\documentclass{amsart}

\usepackage{epsfig}
\usepackage{diagrams}
\newdiagramgrid{small}{1.5,1.5}{0.5,0.5}

\theoremstyle{definition}
\newtheorem{defin}{Definition}[section]
\theoremstyle{remark}
\newtheorem{quest}[defin]{Question}
\newtheorem{prob}[defin]{Problem}

\newtheorem{ex}[defin]{Example}
\newtheorem{rem}[defin]{Remark}
\theoremstyle{plain}
\newtheorem{cor}[defin]{Corollary}
\newtheorem{thm}[defin]{Theorem}
\newtheorem{lem}[defin]{Lemma}
\newtheorem{prop}[defin]{Proposition}
\newtheorem{notat}[defin]{Notation}

\newcommand{\theoremx}[1]{\begin{thm} #1 \end{thm}}

\newcommand{\propositionx}[1]{\begin{prop} #1 \end{prop}}

\newcommand{\corollaryx}[1]{\begin{cor} #1 \end{cor}}

\newcommand{\lemmax}[1]{\begin{lem} #1 \end{lem}}

\newcommand{\definitionx}[1]{\begin{defin} #1 \end{defin}}
\newcommand{\questionx}[1]{\begin{quest} #1 \end{quest}}
\newcommand{\remarkx}[1]{\begin{rem} #1 \end{rem}}

\newcommand{\examplex}[1]{\begin{ex} #1 \end{ex}}
\newcommand{\problemx}[1]{\begin{prob} #1 \end{prob}}
\newcommand{\notationx}[1]{\begin{notat} #1 \end{notat}}

\newcommand{\PL}{\mathrm{PL}}
\newcommand{\Fix}{\mathrm{Fix}}

\title[The Simultaneous Conjugacy Problem in Groups of PL-functions]
{The Simultaneous Conjugacy Problem in Groups of Piecewise Linear Functions}

\thanks{$^\ast$The first author was partially supported by  grants
from National Science Foundation DMS~0600244, 0635607 and~0900932.}

\author{Martin Kassabov$^\ast$ \and Francesco Matucci}


\dedicatory{{\small $^\ast$\emph{Cornell University, Department of Mathematics,
Malott Hall, Ithaca, NY 14853, USA} \\
{\tt{kassabov@math.cornell.edu}}
\\ \bigskip
$^\dagger$\emph{University of Virginia, Department of Mathematics,
Kerchof Hall, Charlottesville, VA 22904, USA} \\
{\tt{fm6w@virginia.edu}}
\\ 
\bigskip
\bigskip
This work is dedicated to the memory of Will Davis}}

\begin{document}

\begin{abstract}
Guba and Sapir asked 
if
the simultaneous conjugacy problem was solvable in Diagram Groups or,
at least, for Thompson's group $F$. 
We give a
solution to the latter question using elementary techniques which 
rely purely on the description
of $F$ as the group of piecewise linear
orientation-preserving homeomorphisms of the unit interval.
The techniques we develop extend the ones used by Brin and Squier allowing us
to compute roots and centralizers as well. Moreover, these techniques
can be generalized to solve the same question in larger groups of
piecewise-linear homeomorphisms.
\end{abstract}

\maketitle

\section{Introduction}

Richard Thompson's group $F$ can be defined by the following presentation:
\[
F = \langle x_0, x_1, x_2, \ldots \mid x_n x_k = x_k x_{n+1}, \forall \, k<n\rangle
\]
This group was introduced and studied by Thompson in the
1960s. The standard
introduction to $F$ is~\cite{cfp}. The group $F$ can be regarded as a subgroup of the group of
piecewise linear self-homeomorphisms of the unit interval and this is the point of view
that we will adopt throughout the paper, and that we will introduce in detail in
Section~\ref{sec:strategy}.

We say that a group $G$ has \emph{solvable ordinary conjugacy problem}
if there is an algorithm such that, given any two elements $y,z \in G$,
we can determine whether there is,
or not, a $g \in G$ such that $g^{-1}yg=z$. Similarly, for a
fixed $k \in \mathbb{N}$, we say that the group $G$ has \emph{solvable $k$-simultaneous
conjugacy problem} if there is an algorithm such that, given any two $k$-tuples of elements in $G$,
$(y_1,\ldots,y_k),(z_1,\ldots,z_k)$, can determine whether there is, or not, a $g \in G$ such that
$g^{-1} y_i g = z_i$ for all $i=1, \ldots, k$. For both these problems, we say that there is an
\emph{effective solution} if the algorithm produces such an element $g$,
in addition to proving its existence.

This problem was studied before for various classes of groups.
The 
simultaneous conjugacy problem has been proved to be solvable for the matrix
groups $\mathrm{GL}_n(\mathbb{Z})$ and $\mathrm{SL}_n(\mathbb{Z})$
by Sarkisyan in 1979 in~\cite{sarkisyan} and independently by
Grunewald and Segal in 1980 in~\cite{grunsegal}.
In 1984 Scott constructed examples of finitely presented 
groups that have an unsolvable conjugacy problem in~\cite{scottelizabeth}.
In 1976 Collins showed in~\cite{collins} that 
the solvability of the conjugacy problem does not imply the solvability of simultaneous conjugacy problem.
More recently, in their 2005 paper~\cite{bridhowie},
Bridson and Howie constructed examples of finitely presented groups
where the ordinary conjugacy problem is solvable, but the $k$-simultaneous
conjugacy problem is unsolvable for every $k \ge 2$.


The ordinary conjugacy problem for $F$
was addressed by Guba and Sapir~\cite{gusa1} in 1997, who
reduced the solution of the conjugacy problem for diagram groups to the solution of the word problem
in the corresponding semigroup, solving this last problem for $F$ and many similar groups.
Their solution, for general diagram groups,
reduces the problem to the isomorphism problem
of planar graphs.
We mention here relevant related work: in 2001 Brin and Squier in~\cite{brin2}
produced a criterion for describing conjugacy classes in $\PL_+(I)$,
the group of all piecewise-linear orientation preserving self-homeomorphisms
of the unit interval with only finitely many breakpoints,
that contains $F$ as a proper subgroup.
In 2007 Gill and Short~\cite{ghisho1} extended this criterion to work in $F$,
thus finding another way to characterize conjugacy classes
from a piecewise linear point of view. Using an approach similar to
Guba and Sapir's original solution, in 2007 Belk and Matucci~\cite{bema1}
produced a unified solution of the conjugacy problem for all Thompson's groups $F,T$ and $V$.

In 1999, Guba and Sapir~\cite{gusa2} posed the question
of whether or not the simultaneous conjugacy problem was solvable for diagram groups.
Some of the results of the present paper are already known, but we deduce all of them
using our tools. We will show that our techniques can be used on a large class of
groups of piecewise linear homeomorphisms. 

\theoremx{
\label{A}
Thompson's group $F$ has a solvable
$k$-simultaneous conjugacy problem, for every $k \in \mathbb{N}$.
There is an algorithm which produces an effective solution and
enumerates all possible conjugators.}
The same algorithm also solves the $k$-simultaneous
conjugacy problem in many ``Thompson-like'' subgroups of $\PL_+(I)$ (see
Subsection~\ref{ssec:notation} for the precise definition).

As an application of the proof of the above
Theorem we have the following corollaries:

\theoremx{
\label{B}
For an element $x \in F$, we denote by $C_F(x)$ the centralizer of $x$ in $F$. Then:
\begin{enumerate}
\item[(i)] $C_F(x) \cong  F^m \times \mathbb{Z}^n$, for some numbers $0 \le m \le n + 1$.
\item[(ii)] $x\in F$ has a finite number of roots, which can be effectively computed.
\item[(iii)] The centralizer of any finitely generated subgroup $A\subset F$ decomposes as the direct
product of the groups $C_i$, where each $C_i$ is either trivial, infinite cyclic or isomorphic to $F$.
\item[(iv)] The intersection of any number $k \ge 2$ centralizers of elements of $F$ is equal to the
intersection of 2 centralizers.
\end{enumerate}}

Parts of the previous theorem were already proved either in the setting of $F$ or in that of $\PL_+(I)$:
in particular parts (i) and (ii) were proved by Guba and Sapir in \cite{gusa1} for $F$ 
and by Brin and Squier for $\PL_+(I)$ in \cite{brin2}.
All the previous results can be suitably rephrased 
for a large class of subgroups of $\PL_+(I)$
(see Subsection~\ref{ssec:notation} for the precise definition).

\medskip

The
paper is organized as follows: in Section~\ref{sec:strategy} we will define the groups $\PL_{S,G}(I)$
that generalize Thompson's group $F$ and give an outline of the
solution of simultaneous conjugacy problem.
In Section~\ref{sec:main-machine} we introduce the main
algorithm to create candidate conjugators. In Section~\ref{sec:cenPL}
we compute centralizers and roots.
In Section~\ref{ssec:diagonls} we show how to build an approximate conjugator which makes
the fixed point set of $y$ and $z$ coincide.
In Section~\ref{sec:applications} we get the solution of the ordinary conjugacy
problem and a variation of it, the \emph{power conjugacy problem}.
In Section~\ref{sec:simconj} we describe how to reduce the simultaneous
conjugacy problem to a special instance of the ordinary conjugacy problem.
In Section~\ref{sec:examples} we show interesting instances
where the simultaneous conjugacy problem can be solved.

\section{The idea of the argument}\label{sec:strategy}

In this section we describe the groups that we will study and
outline the steps of our proof. This section
is intended to give a quick overview of the results that we will prove in the later sections.

\subsection{Notations\label{ssec:notation}}

We introduce here the notation that will be used across the paper.
Let $I=[0,1]$ be the unit interval. We define $\PL_+(I)$ to be the
group of piecewise linear%
\footnote{By piecewise linear we mean piecewise affine, although this abuse of language is now common.}
orientation-preserving homeomorphisms of unit interval
into itself,
with finitely many breakpoints of the derivative function such that slopes are positive real numbers.
The product of two elements is given by the composition of functions.

\medskip

One can impose additional the requirements on the breakpoints
and the slopes to define subgroups of $\PL_+(I)$.
Let $S$ be an additive subgroup of $\mathbb{R}$ containing $1$, 
let $U(S)$ denote
the multiplicative group 
$\{g \in \mathbb{R}^* | gS = S \mbox{ and } g > 0\}$ 
and let $G$ be a subgroup of
$U(S)$. Thus, $S$ is a module over the ring $\mathbb{Z}[G]$.
We define $\PL_{S,G}(I)$ to be the
subgroup of $\PL_+(I)$
consisting of all functions $f$  such that
the breakpoints are in the subgroup $S$ and
the slopes are in the subgroup $G$. We observe that 
if the group $G$ is trivial, then so is $\PL_{S,G}(I)$. Therefore in the rest of the
paper we will assume that $G$ is nontrivial, which implies that $S$
is dense in $\mathbb{R}$ (with respect to the usual topology).

If $G=U(S)$ we
write $\PL_S(I)$, instead of $\PL_{S,G}(I)$.
If $S=\mathbb{R}$, then $\PL_S(I)=\PL_+(I)$.
For the special case $S=\mathbb{Z}\left[\frac{1}{2}\right]$,
we denote the group $\PL_{\mathbb{Z}{\left[\frac{1}{2}\right]}}(I)$
by $\PL_2(I)$. The group $\PL_2(I)$ is also known as
\emph{Thompson's group $F$} and it is isomorphic to the group $F$
defined in the introduction (see~\cite{cfp} for a proof).
\footnote{The family of groups $\PL_{S,G}(I)$ was first introduced by Bieri and Strebel in
~\cite{bieristrebel1} and was later popularized through the work of Stein ~\cite{stein1}.}
We remark that in order to make some calculations possible
inside the module $S$ and its quotients, we need to ask for some requirements
to be satisfied by $S$ from the computability standpoint
(like the existence of black box algorithms
for performing the basic operations in $S$).
These will be explicitly stated in section~\ref{sec:comp_req}
and will be assumed throughout this paper.

To attack the ordinary and the simultaneous conjugacy problems,
we will split the study into that of some families of 
functions inside $\PL_+(I)$.
The reduction to these subfamilies will come from the study
of the fixed point subset of the interval $I$ for a function $f$.

\remarkx{
\label{thm:define-PLSG}
We would like to define the group $\PL_{S,G}(J)$,
where $J=[\eta,\zeta]$ is
any interval contained in $I$. We consider
the  group of restrictions of functions in $\PL_{S,G}(I)$
fixing the endpoints of $J$:
\[
\PL_{S,G}^{\mathrm{Rest}}(J):=
\left\{f|_J  \mid  f \in \PL_{S,G}(I), f(\eta)=\eta, f(\zeta)=\zeta \right\}
\]
In general, it is not true that $\PL_{S,G}^{\mathrm{Rest}}(J)$
is a subgroup of $\PL_{S,G}(I)$. Moreover, there is no natural embedding of $\PL_{S,G}^{\mathrm{Rest}}(J)$
into $\PL_{S,G}(I)$ such that the restriction of the image of a function is the initial function
(see also Remark~\ref{thm:multiple-orbits-of-points}).
If the endpoints of $J$ are in $S$, we will denote the group
$\PL_{S,G}^{\mathrm{Rest}}(J)$ with $\PL_{S,G}(J)$.%
\footnote{There is another natural way to define $\PL_{S,G}(J)$: consider
the subgroup of functions of $\PL_{S,G}(I)$
which fix the endpoints of $J$ and are the identity on $I \setminus J$:
\[
\PL_{S,G}^{\Fix=I \setminus J}(J):=
\left\{f \in \PL_{S,G}(I) \mid f(t)=t, \forall t \in I \setminus J\right\}.
\]
We observe that by definition $\PL_{S,G}^{\Fix=I \setminus J}(J)$
is a subgroup of $\PL_{S,G}(I)$.
In the case where the endpoints of $J$ are contained in $S$,
the two definitions coincide, i.e.,
$\PL_{S,G}^{\Fix = I \setminus J}(J) \cong \PL_{S,G}^{\mathrm{Rest}}(J)$,
and thus the group $\PL_{S,G}^{\mathrm{Rest}}(J)$
can be regarded as a subgroup of $\PL_{S,G}(I)$.
}
}

\remarkx{
\label{thm:endpoints-in-S}
Throughout the paper we will always assume the interval $J$
to have endpoints in $S$ (with the only exception of Lemma \ref{thm:nec-cond-trans}).
For the special case $S = \mathbb{Z}{\left[\frac{1}{2}\right]}$,
it is straightforward to verify that
$\PL_2(J) \cong \PL_2(I)$. We observe that the analogous fact may
not be true for the groups $\PL_{S,G}(I)$
(see Remark~\ref{thm:multiple-orbits-of-points}).}

For a function $f \in \PL_{S,G}(J)$ we define the fixed point set
on the interval $J$ as
$$
\Fix_J(f):=\{t \in J \;\big\vert\; f(t)=t \},
$$
which is a closed set. It follows from the definition that
$\Fix_J(f)$ is a union of finitely many intervals with end points in $S$
and finitely many ``isolated'' points.
We will often simplify the notation by dropping the subscript $J$.
The motivation for introducing this subset is easily explained ---
if $y,z \in \PL_+(J)$ are conjugate through
$g \in \PL_+(I)$ and $t \in (\eta,\zeta)$ is such that $y(t)=t$ then
$z(g^{-1}(t))=(g^{-1}yg)(g^{-1}(t))=g^{-1}(t)$, that is,
if $y$ has a fixed point then $z$ must have a fixed point.

\begin{defin}
We define $\PL_{S,G}^<(J)$ and $\PL_{S,G}^>(J)$ to be the set of all functions
in $\PL_{S,G}(J)$ with graph below the diagonal, respectively above the diagonal.
Following Brin and Squier~\cite{brin2}, we define
a function in $x \in \PL_{S,G}(J)$ to
be a \emph{one-bump function} if either
$x \in \PL_{S,G}^>(J)$ or $x \in \PL_{S,G}^<(J)$.
\end{defin}

\begin{figure}
 \begin{center}
  \includegraphics[height=5cm]{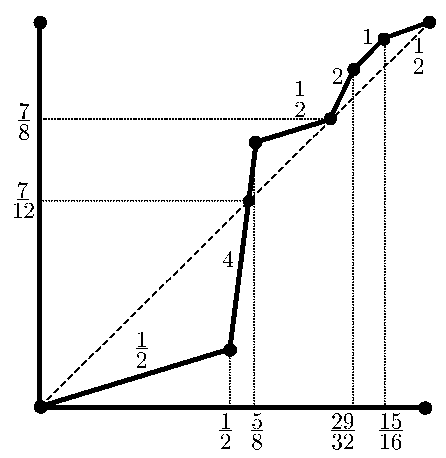}
 \end{center}
 \caption{A function in $\PL_2(I)$ with a non-dyadic fixed point.}
 \label{fig:rational-intersection}
\end{figure}

In general it is not true that  if $f \in \PL_{S,G}(I)$ then
$\Fix(f)\subseteq S$, but 
$\Fix(f)$ is always a
subset of the ``field of fractions'' of $S$.
The example in figure~\ref{fig:rational-intersection}
shows a function in $\PL_2(I)$
with a non-dyadic rational fixed point.
In order to avoid working in intervals $J$ where the endpoints may not be in $S$,
we introduce a new definition of boundary which deals with this situation:
for a subset $X \subseteq [0,1]$,
we define
$$
\partial_S X:=\partial X \cap S
$$
%
where $\partial X$ denotes the usual topological boundary of
$X$ inside $\mathbb{R}$.
For the special case $S=\mathbb{Z}\left[\frac{1}{2}\right]$
we write $\partial_2 X$.
We are going to apply this definition to the set
$X=\Fix(f)$ so that $\partial \Fix(f)$ and
$\partial_S \Fix(f)$
will always be finite.

\definitionx{\label{thm:define-PLSG-no-intersection}
The $\PL_{S,G}^0(J) \subseteq \PL_{S,G}(J)$ will denote
set of functions $f \in \PL_{S,G}(J)$ such that
the set $\Fix(f)$ does not contain elements of $S$
other than the endpoints of $J$,
i.e., $\Fix(f)$ is discrete and
$\partial_S \Fix(f)=\partial_S J=\partial J$.
The elements in $\PL_{S,G}^0(J)$ will be called 
\emph{almost one-bump function}, although their graphs 
have several bumps in general.}

\subsection{Outline of the strategy}
\label{ssec:outline-strategy}

We are now going to describe the general steps and reductions
of the algorithm to solve the simultaneous conjugacy problem in
the groups $\PL_{S,G}(I)$. Most of the time we will work in the larger group
$\PL_+(I)$ and we will then say what is necessary to generalize
the argument to $\PL_{S,G}(I)$. The following outline describes the
correct order of the steps needed to solve the problem,
however we will start Section~\ref{sec:main-machine}
by describing the central tool of the paper (the ``Stair Algorithm'')
which is used in Step 2. Let $x,y,z \in \PL_{S,G}(I)$.
\medskip

\noindent \textbf{Step 1.}
\emph{Find a $g \in \PL_{S,G}(I)$ such that $\Fix(y)=g(\Fix(z))$.}
The set $\Fix(x)$
consists of a disjoint union of a finite number of closed intervals and
isolated points, because every $x \in \PL_{S,G}(I)$ has only finitely
many breakpoints.
As mentioned before, if $g^{-1}yg=z$, then $\Fix(y)=g(\Fix(z))$.
Thus, as a first step
we need to know if, given $y$ and $z$, there exists a $g \in \PL_{S,G}(I)$
such that $\Fix(g^{-1}yg) = g(\Fix(y))=\Fix(z)$.
In Section~\ref{ssec:diagonls} we show an algorithm which determines
whether or not there exists a ``candidate''
conjugator $g_\ast$ such that $\Fix(g_\ast^{-1} y g_\ast)= \Fix(z)$.
We then study the conjugacy problem problem for
$g_\ast^{-1} y g_\ast$ and $z$.

\medskip

\noindent  \textbf{Step 2.}
\emph{Solve the conjugacy problem if $\Fix(y)=\Fix(z)$}.
In this case
$\partial_S \Fix(y) = \partial_S \Fix(y)=
\{\alpha_1,\ldots,\alpha_n\}$. It is easy to see that any conjugator
fixes the points $\alpha_i$, this we need to
look for conjugators in $\PL_{S,G}([\alpha_i,\alpha_{i+1}])$
of the restrictions of $y$ and $z$ to $[\alpha_i,\alpha_{i+1}]$.
Thus, we can reduce the conjugacy problem to the intervals where
$y$ and $z$ are ``almost one-bump functions'', more precisely
they are either in $\PL^0_{S,G}(J)$ or equal to the identity.
The case $y=z=id$ is trivial, this we can assume that
 $y$ and $z$ are ``almost one-bump'' functions%
\footnote{One needs to be a bit more careful since $y$ and $z$
can have fixed points in the interval which do not lie in $S$.}.
This case will be dealt with through a procedure called the
``stair algorithm'' that we provide in Subsection~\ref{ssec:stair}.

\medskip

\noindent \textbf{Step 3.}
\emph{Describe the intersection of centralizers of elements and derive
a solution to the conjugacy problem}.
Finding centralizers $g$ of an element $y$ is equivalent to
find all elements $g$ such that $g^{-1}yg=y$.
Using similar techniques we can also classify the structure of
intersection of centralizers, which will be useful for the last step.
Since the set of all conjugators for $y$ and $z$ is given
by a particular conjugator times an element in the centralizer of $y$,
Steps 1, 2 and 3 give us a solution to the conjugacy problem.

\medskip

\noindent \textbf{Step 4.}
\emph{Reduce the simultaneous conjugacy problem to a
``restricted'' conjugacy problem.}
It can be seen that the simultaneous conjugacy problem is
equivalent to solving the conjugacy problem
for two elements $y$ and $z$ with the restriction that the conjugator
$g$ must lie in the intersection
of centralizers of some elements $x_1, \ldots, x_k$.
In Section~\ref{sec:simconj} we will show how to build
such a conjugator, if it exists, following the previous steps.

\section{Computational Requirements}
\label{sec:comp_req}

In order to effectively solve to conjugacy and the simultaneous conjugacy
problem in the groups $\PL_{S,G}$ we need to assume that the additive group
$S$ and the multiplicative group $G$ satisfy some computational requirements.
First, we will assume that there is some representation of the elements in $S$ and $G$
in some data structure $M$.%
\footnote{
Usually the elements are represented by some finite strings over a given alphabet,
if this is the case we will require that the sets $S$ and $G$ are countable.
However our algorithms does not depend on the data structure $M$.}
Then we need to be able to perform the basic operations in $S$ and $G$, thus
we require that we are give some ``oracles'' which perform the
following operations:
\begin{itemize}
\item determine if $m\in  M$ represents an element in $S$ and/or $G$;
\item determine if $m,m'\in M$ represent the same element in $S$ and/or $G$;
\item perform the basic operations (additions and substraction)
in $S$;
\item given two elements in $S$, determine which one is bigger;
\item given an element of $S$ and a rational number, determine which one is bigger;
\item construct an element in $S$ in any given non-empty
open interval;
\item perform the basic operations (multiplication, division) in $G$;
\item perform multiplication between the elements in $G$ and $S$.
\end{itemize}
Using these oracles, one can construct a data structure which represents the
elements in the group $\PL_{S,G}$ and new oracles which perform the group
operations.

The following additional oracles are needed for the algorithms describes in
section~\ref{ssec:diagonls} (here $\mathcal{I}$ denote the subgroup of
$S$ generated by $(g-1)s$ for $s\in S$ and $g\in G$):
\begin{itemize}
\item given $g\in G$ and $s\in S$ determine if $s/(g-1)$ is an element in $S$;
\item effective solution of the membership problem in the submodule $\mathcal{I}$,
i.e., given $s\in S$, an oracle determines if $s\in I$ or not and
if $s\in \mathcal{I}$ it produces elements $s_i$ and $g_i$
such that $s = \sum (g_i-1) s_i$;
\item effective solution of the congruence 
$sG=s'G \pmod{(t-1) \mathcal{I}}$, i.e.,
given $s,s'$ and $t$ an oracle constructs a solution of the congruence
or determines that it has no solutions.
\end{itemize}
These oracles allow us to solve effectively solve the conjugacy problem
in the group $\PL_{S,G}(J)$, but for an effective solution of the
simultaneous conjugacy problem we need another oracle
\begin{itemize}
\item effective solution of the equation  $a^k = b c^i$, where $k,i\in \mathbb{Z}$, 
i.e.,  given $a,b,c \in G$ construct an integer  solution of the equation
or determine that there cannot be any.
\end{itemize}

\section{The Stair Algorithm}
\label{sec:main-machine}

In this Section we carry out the second step of the
strategy described in Subsection~\ref{ssec:outline-strategy}
by restricting our study to a square where the given functions have
``no relevant'' intersection with the diagonal, and showing
how to build possible candidates for conjugator. Our goal for this
section is, essentially, to solve the conjugacy problem in $\PL_+^<(J)$
where we do not pay attention to the intersection with the diagonal.
Our methods extend the results of Brin and Squier~\cite{brin2},
who develop a technique similar to our algorithm.
In this section we develop an algorithm, which is allows us to recover
Brin and Squier's analysis and extend it to the case of $\PL_{S,G}(J)$,
together with a description of the intersection of centralizers.

\subsection{The Linearity Boxes}
\label{ssec:linearitybox}

This subsection and the following one will deal with functions in
$\PL_+(J)$ for an interval $J=[\eta,\zeta]$:
we will reuse them in the discussion on $\PL_{S,G}(I)$.
We start by making the following observation:
the map $\PL_+(J) \to \mathbb{R}_+$ which sends a function $f$ to $f'(\eta^+)$
is a group homomorphism. The very first thing to check,
if $y$ and $z$ are to be conjugate through a $g \in \PL_+(J)$
in neighborhoods of the endpoints of $J$,
the following trivial lemma says that this can happen only if
the graphs of $y$ and $z$ coincide near the end points of $J$.

\begin{figure}
 \begin{center}
  \includegraphics[height=5cm]{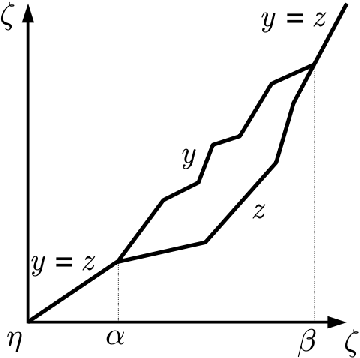}
 \end{center}
 \caption{$y$ and $z$ coincide around the endpoints.}
 \label{fig:startendequal}
\end{figure}

\lemmax{
\label{thm:startendequal}
Given three functions $y,z,g \in \PL_+(J)$ such that $g^{-1} y g= z$,
there exist $\alpha, \beta \in (\eta,\zeta)$ such that $y(t)=z(t)$, for all
$t \in [\eta,\alpha] \cup [\beta,\zeta]$ (refer to figure~\ref{fig:startendequal}).
}

The next lemma gives us that any function $g\in \PL_+(J)$ which conjugates
$y$ to $z$ needs to be linear in a specific neighborhood of each endpoint of $J$,
which depends only on $y$ and $z$. This lemmas is main ingredient which allows us
to extend the methods of Brin and Squier~\cite{brin2} to get constructive
solution of the conjugacy problem.

\begin{figure}
 \begin{center}
  \includegraphics[height=5cm]{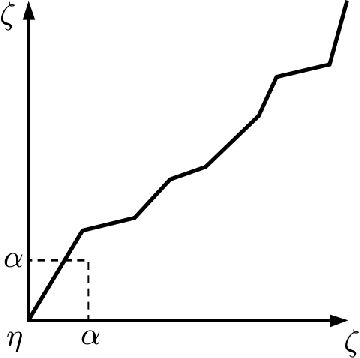}
 \end{center}
 \caption{Initial linearity box.}
 \label{fig:linearityboxes}
\end{figure}

\lemmax{
\label{thm:initial-box}
Suppose $y,z,g \in \PL_+(J)$ and $g^{-1}yg=z$.
Let $\varepsilon>0$ and $y'(\eta^+)=z'(\eta^+)=c > 1$ satisfy
$$
y(t) - \eta = z(t)-\eta = c (t - \eta) \,\,
\mbox{for } t\in[\eta,\eta+\varepsilon].
$$
Then the graph of $g$ is linear inside the square
$[\eta,\eta+\varepsilon]\times[\eta,\eta+\varepsilon]$ 
(see figure~\ref{fig:linearityboxes}).
}
\begin{proof}
We can rewrite the conclusion of this lemma, by saying that, if we define
$$
\tilde\varepsilon = \sup \{r \mid g \mbox{ is linear on } [\eta,\eta+r]\},
$$
then $\eta +\tilde\varepsilon\ge \min \{g^{-1}(\eta+\varepsilon),\eta+\varepsilon\}$.
Assume the contrary, let
$\tilde \varepsilon < \varepsilon$ and
$\eta + \tilde\varepsilon < g^{-1}(\eta + \varepsilon)$ and write
$g(t)- \eta =\gamma (t-\eta)$
for  $t\in [\eta,\eta+\varepsilon]$, for some constant $\gamma > 0$.
Let $0 \le \sigma <1$ be any number,
since  $\tilde \varepsilon < \varepsilon$,
we have $\eta+\sigma \tilde\varepsilon<\eta + \varepsilon$ and
so $y$ is linear around $\eta+\sigma \tilde\varepsilon$:
\[
g(y(\eta+ \sigma \tilde \varepsilon))=g(\eta+c \sigma \tilde \varepsilon).
\]
On the other hand, since $\eta + \tilde \varepsilon < g^{-1}(\eta + \varepsilon)$,
it follows that
$g(\eta+\sigma\tilde \varepsilon) < g(\eta+\tilde\varepsilon)< \eta+\varepsilon$
and so $z$ is linear around the point $g(\eta+ \sigma \tilde\varepsilon)=
\eta +\gamma \sigma \tilde\varepsilon$:
$$
z(g(\eta+ \sigma \tilde\varepsilon))= z(\eta +\gamma \sigma \tilde\varepsilon)=
\eta+ c \gamma \sigma \tilde\varepsilon.
$$
Since $gy=zg$, we can equate the previous two equations and write
$g(\eta+c \sigma \tilde\varepsilon)=\eta+\gamma c \sigma \tilde\varepsilon$,
for any number $0 \le \sigma <1$.
If we choose $1/c<\sigma<1$, we see that $g$ must be linear on the interval
$[0, c \sigma \tilde\varepsilon]$,
where  $c \sigma \tilde\varepsilon> \tilde\varepsilon$.
This is a contradiction to the definition of $\tilde\varepsilon$.
\end{proof}

We observe that the Lemma also holds when $z'(\eta^+)=y'(\eta^+)=c<1$
by applying it to the homeomorphisms $y^{-1},z^{-1}$.
Thus we can replace the condition $z'(\eta^+)=y'(\eta^+)=c>1$ with
$z'(\eta^+)=y'(\eta^+)\not =1$.
Note that Lemma~\ref{thm:initial-box} has an analogue
for the points close to other endpoint of $J$:

\remarkx{
\label{thm:final-box}
Let $y,z,g \in \PL_+(J)$. Suppose $g^{-1}yg=y$,
if there exist $\beta \in J$ and $c<1$ such that
$y(t)=z(t)=c \cdot (t-\zeta) + \zeta$
on $[\beta,\zeta]$,
then the graph of $g$ is linear inside the square
$[\beta,\zeta]\times[\beta,\zeta]$.
}

Lemma~\ref{thm:initial-box} does not hold
when the initial slope of $y$ and $z$ are equal to $1$,
because any function $g$ with a support sufficiently close to the
end points will conjugate $y$ to itself.

\subsection{The Stair Algorithm for $\PL_+^<(J)$}
\label{ssec:stair}


This subsection will deal with the main construction of this paper.
We show how, under certain hypotheses, if there is a conjugator, then it is unique.
On the other hand, we give a construction of such a conjugator, if it exists.
Given two elements $y,z$ the set of their conjugators is
a coset of the centralizer of one of them, thus it makes sense
to start by deriving properties of centralizers.

The first several Lemmas show that if $y$ and $z$ are one bump functions
in $\PL_+(J)$ then the graphs of the conjugators do not intersect.

\lemmax{
Let $z \in \PL_+(J)$. Suppose there exist $\eta \le \lambda \le \mu \le \zeta$
such that $z(t) \le \lambda$, for every $t \in [\eta,\mu]$. Suppose further that
$g \in \PL_+(J)$ is such that
$g(t)=t$, for all $t \in [\eta,\lambda]$ and
$g^{-1}z g(t)=z(t)$, for all $t\in [\eta,\mu]$.
Then $g(t)=t$, for all $t \in [\eta,\mu]$.
}

\begin{proof}
The equation $g^{-1}z g(t)=z(t)$ implies that
$g(t)=z^{-1} g z(t)$ for all $t\in [\eta,\mu]$.
Since $z(t) \le \lambda$ and $g(x)=x$ for all $x\leq \lambda$, we have
$$
g(t)=z^{-1}( g (z(t))) = z^{-1}(z(t)) = t. \qedhere
$$
\end{proof}

\corollaryx{
\label{thm:unique-partial}
Let $z \in \PL_+^<(J)$ and $g \in \PL_+(J)$ be such that
$g'(\eta^+)=1$,
and
$g^{-1}z g=z$.
Then $g(t)=id$, the identity map.
}

\begin{proof}
Since $g'(\eta^+)=1$, we have $g(t)=t$ for all $t\in [\eta, \eta + \varepsilon]$.
Applying the previous lemma several times gives that $g(t)=t$ for all
$t \in [\eta, z^{-k}(\eta + \varepsilon)]$. Since $z \in \PL_+^<(J)$ we have that
$\lim_{k\to\infty }z^{-k}(\eta + \varepsilon) = \zeta$, therefore $g(t) = t$
for $t\in J$.
\end{proof}

\lemmax{
\label{sec:unque}
Let $z \in \PL^<_+(J)$.
Let $C_{\PL_+(J)}(z)$ be the centralizer of $z$ in $\PL_+(J)$.
Then $\varphi_z$, defined below, is an injective group homomorphism.
$$
\begin{array}{lrcl}
\varphi_z: & C_{\PL_+(J)}(z) & \longrightarrow & \mathbb{R_+} \\
           & g              & \longmapsto     & g'(\eta^+).
\end{array}
$$
}

\begin{proof}
Clearly $\varphi_z$ is a group homomorphism.
Suppose that there exists two elements $g_1,g_2 \in C_{\PL_+(J)}(z)$ such that
$\varphi_z(g_1)=\varphi_z(g_2)$, then $g^{-1}_1 g_2$ has a slope $1$ near $\eta$
and by the previous Lemma is equal to the identity. Therefore $g_1=g_2$,
which proves the injectivity of $\varphi_z$.
\end{proof}


\lemmax{
\label{cong:unique}
Let $y,z \in \PL^<_+(J)$,
let $C_{\PL_+(J)}(y,z)=\{g \in \PL_+(J) \; | \; y^g=z\}$ be the set
of all conjugators and let $\lambda$ be an interior point of $J$.
Then the two maps $\varphi_{y,z}$ and $\psi_{y,z,\lambda}$
satisfy
$$
\begin{array}{lrcl}
\varphi_{y,z}: & C_{\PL_+(J)}(y,z) & \longrightarrow & \mathbb{R_+} \\
               & g                & \longmapsto     & g'(\eta^+)
\\
\psi_{y,z,\lambda}:      & C_{\PL_+(J)}(y,z)   & \longrightarrow & J \\
             & g                & \longmapsto     & g(\lambda).
\end{array}
$$
\begin{itemize}
\item[(i)] $\varphi_{y,z}$ is an injective map.

\item[(ii)] There is a map $\rho_\lambda:J \to \mathbb{R}_+$
such that the following diagram commutes:
\begin{displaymath}
\begin{diagram}[grid=small]
                  &                            & J \\
                  & \ruTo^{\psi_{y,z,\lambda}} & \dTo_{\rho_\lambda} \\
C_{\PL_+(J)}(y,z) & \rTo_{\varphi_{y,z}}       & \mathbb{R}_+ \\
\end{diagram}
\end{displaymath}

\item[(iii)]$\psi_{y,z,\lambda}$ is injective.
\end{itemize}
}

\begin{proof}
(i) is an immediate corollary of Lemma~\ref{sec:unque}.
(ii) Without loss of generality we can assume that the initial slopes
of $y,z$ are the same (otherwise the set
$C_{\PL_+(J)}(y,z)$ is obviously empty and any map will do).
We define the map $\rho_\lambda:J \to \mathbb{R}_+$ as
\[
\rho_\lambda(\mu)= \lim_{n \to \infty} \frac{y^n(\mu)-\eta}{z^n(\lambda)-\eta}
\]
The above limit exists, because the sequence stabilizes under the assumptions
$y,z \in \PL^<_+(J)$ and $y'(\eta)=z'(\eta)$.

To prove that the diagram commutes we define $\mu=g(\lambda)$ and observe that
$y^{n}(\mu)\underset{n \to \infty}{\longrightarrow}\eta$ and
$z^{n}(\lambda)\underset{n \to \infty}{\longrightarrow}\eta$.
By hypothesis $y(\mu)=g(z(\lambda))$ so that $g(z^n(\lambda))=y^n(\mu)$,
for every $n \in \mathbb{Z}$.
Since $g$ fixes $\eta$ we have
$$
g(t)= g'(\eta^+)(t-\eta)+\eta
\mbox{ on a small interval } [\eta,\eta + \varepsilon],
$$
where $\varepsilon$ depends only on $g$. Let $N=N(g) \in \mathbb{N}$
be large enough, so that the numbers
$y^N(\lambda), z^N(\lambda) \in (\eta,\eta+\varepsilon)$.
This implies that, for any $n\ge N$
$$
y^n(\mu)=g(z^n(\lambda))=g'(\eta^+)(z^n(\lambda)-\eta)+\eta
$$
and so then
\[
\varphi_{y,z}(g)=g'(\eta^+)=
\frac{y^n(\mu) -\eta}{z^n(\lambda)-\eta}=\rho_{\lambda}(\psi_{y,z,\lambda}(g)).
\]
\noindent (iii) Since $\varphi_{y,z}=\rho_{\lambda}\psi_{y,z,\lambda}$
is injective by part (i), then $\psi_{y,z,\lambda}$ is also injective.
\end{proof}

\begin{figure}
 \begin{center}
  \includegraphics[height=5cm]{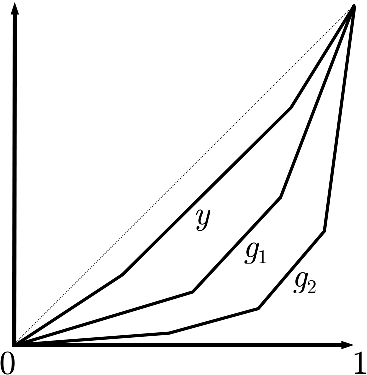}
 \end{center}
 \caption{Two elements $g_1,g_2$ centralizing a map $y \in \PL^<_+(J)$.}
 \label{fig:centralizer-not-intersecting}
\end{figure}

\remarkx{
Lemma~\ref{cong:unique} shows that for $y \in \PL^<_+(J)$ the
graphs of the elements in the centralizer $C_{\PL^<_+(J)}(y)$ do not
intersect, see figure \ref{fig:centralizer-not-intersecting}.}

The main tool of this subsection is the {\bf Stair Algorithm}.
This procedure builds a conjugator (if it exists) with
a given fixed initial slope. 
In order for $y$ and
$z$ to be conjugate, they must have the same initial slope;
by Lemma~\ref{thm:initial-box} this determines uniquely
the first piece of a possible conjugator given the initial slope.
Then we ``walk up the first step of the stair''
(Lemma~\ref{thm:extend-conjugator}):
we identify $y$ and $z$ inside a rectangle next to the linearity box,
by taking a suitable conjugator. We then repeat and walk up more rectangles,
until we ``reach the door'' (represented by the final linearity box)
and this happens when a rectangle
that we are building crosses the final linearity box.
This algorithm finishes in finitely many steps because the interval
$J=[\eta,\zeta]$ is bounded.
In other words, we will construct a
``section'' for the map $\varphi_{y,z}$ of Lemma~\ref{cong:unique}.
As a consequence we will also build a ``section'' of the map
$\psi_{y,z,\lambda}$ too.

\lemmax{
\label{thm:extend-conjugator}
Let $y,z \in \PL_+^<(J)$ and $g \in \PL_+(J)$ be functions such that $z=y^g$ and
let $\alpha \in (\eta,\zeta)$.
Then the functions $y$, $z$ and the restriction of
$g$ to the interval $(\eta, \alpha)$ uniquely determines
the restriction of $g$ to the interval $(\eta, z^{-1}(\alpha))$.
}
\begin{proof}
We can rewrite the equation $z=y^g$ as $g = y^{-1} g z$. The value of the right side
of this equation at points inside the interval $(\eta, z^{-1}(\alpha))$ depends only on
$y$, $z$ and restriction of $g$ to the interval $(\eta, \alpha)$. Therefore they
determine uniquely the restriction of $g$ to the interval $(\eta, z^{-1}(\alpha))$.
\end{proof}

\propositionx{
\label{thm:stair-algorithm}
Let $y,z \in \PL_+^<(J)$ and $g \in \PL_+(J)$ be functions such that $z=y^g$.
Then the conjugator $g$ is uniquely determined by its initial slope $g'(\eta)$.
}
\begin{proof}
By Lemma~\ref{thm:initial-box} the graph of the conjugator $g$ is linear in the
box $[\eta,\eta+\varepsilon] \times [\eta,\eta+\varepsilon]$ therefore the slope of
$g'(\eta)$ uniquely determines the restriction of $g$ on the  interval
$(\eta, \alpha)$, for some $\alpha \le \zeta$.
Applying the previous lemma several times we see that this also determines the
restriction of $g$ to the interval $(\eta, z^{-n}(\alpha))$
for any integer $n \geq 0$.
However the function $z$ is in $\PL_+^<(J)$
therefore $\lim_{n \to \infty} z^{-n}(\alpha) = \zeta$, thus
these restrictions determine the function $g$.
\end{proof}

\begin{rem}
Lemma~\ref{thm:extend-conjugator} also holds for any (even non-piecewise linear)
function on the interval $J$.
The argument in the previous proposition gives that for any piecewise linear
functions $y$ and $z$ in $\PL^+(J)$ and any initial slope
there exists a unique conjugating function $g$
which is linear in a neighborhood of the point $\eta$.
Although this function $g$ is piecewise linear on any interval $(\eta,\alpha)$
for any $\alpha < \zeta$, it may not be linear in a neighborhood
of the point $\zeta$ and may not be
piecewise linear on the entire interval $J$.

Using the final linearity box, it is very easy to algorithmically
determine if this function $g$ is a piecewise linear function.
It suffices to construct the restriction of $g$ to the interval $(\eta,\gamma)$
such that the point $(\gamma,g(\gamma))$ is inside the finial linearity box
$[\beta,\zeta]  \times [\beta,\zeta]$.
By remark~\ref{thm:final-box} if there exists a conjugator
then it has to be linear in this box, thus we can determine the rest of
the graph of $g$ and then we can verify if it is indeed a conjugator.
\end{rem}

\corollaryx{
\label{thm:explicit-conjugator}
Let $y,z \in \PL_+^<(J)$, let $[\eta,\alpha]$ be the initial linearity box
and let $q$ be a positive real number. There is an $r \in \mathbb{N}$
such that the unique candidate conjugator
with initial slope $q<1$ is given by
\[
g(t)=y^{-r}g_0z^r(t) \qquad \forall t \in [\eta,z^{-r}(\alpha)]
\]
and linear otherwise, where $g_0$ is any map in $\PL_+(J)$
which is linear in the initial box and
such that $g_0'(\eta^+)=q$.
}

\lemmax{
\label{thm:unique-roots}
Let $y,z \in \PL_+^<(J)$, $g \in \PL_+(J)$ and $n \in \mathbb{N}$.
Then $g^{-1}y g=z$ if and only if
$g^{-1}y^n g=z^n$. }

\begin{proof}
The ``only if'' part is obvious. The ``if'' part follows from the injectivity
of $\varphi_x$ of Lemma~\ref{sec:unque} since $g^{-1}yg$ and $z$
both centralize the element
$z^n$ and have the same initial slope.
\end{proof}

\corollaryx{
\label{thm:psi-surjective}
Let $y,z \in \PL_+^<(J)$, and let $\lambda$ be in the interior of $J$. The map
$$
\begin{array}{lrcl}
\psi_{y,z,\lambda}:      & C_{\PL_+(J)}(y,z)   & \longrightarrow & J \\
                 & g                  & \longmapsto     & g(\lambda).
\end{array}
$$
admits a \emph{section}, i.e., if $\psi_{y,z,\lambda}(g)=\mu \in J$,
then $g$ is unique and can be constructed.
}

\remarkx{
Suppose $y,z \in \PL_+^<(J) \cup \PL_+^>(J)$, then in order to
be conjugate, they will have to be both in $\PL_+^<(J)$ or both in $\PL_+^>(J)$,
because by Lemma~\ref{thm:startendequal} they will have to coincide
in a small interval $[\eta,\alpha]$. Moreover, $g^{-1}yg=z$ if and only
if $g^{-1}y^{-1}g=z^{-1}$, and so, up to working
with $y^{-1},z^{-1}$,
we may reduce to studying the case where they are both in $\PL_+^<(J)$.}

\remarkx{
The stair algorithm
for $\PL_+^<(J)$ can be reversed. This is to say that, given $q$
a positive real number, we can determine whether or not
there is a conjugator $g$ with final slope $g'(\zeta^-)=q$.
The proof is the same: we simply start building $g$ from the final linearity box.}

\remarkx{\label{thm:closure-inside-PLSG}
We mention here that all results of
subsections~\ref{ssec:linearitybox} and~\ref{ssec:stair} can be extended
to the case of $\PL_{S,G}(J)$. All the statements can be reformulated
and proved by replacing every appearence of $\PL_+(J)$ and $\PL_+^<(J)$
with the symbols  $\PL_{S,G}(J)$ and $\PL_{S,G}^<(J)$.
}

The stair algorithm gives a practical way to find conjugators
if they exist and we have chosen
a possible initial slope.
By analyzing the stair algorithm we can see that,
if two elements are in $\PL_{S,G}^<(J)$
and they are conjugate through an element with
initial slope in $G$ then the conjugator
is an element of $\PL_{S,G}(J)$.

\corollaryx{\label{thm:conjugator-still-in-PLSG}
Let $y,z \in \PL_{S,G}^<(J)$, $g \in \PL_+(J)$
such that $y^g=z$ and $g'(\eta^+)\in G$.
Then $g \in \PL_{S,G}(J)$.}

We conclude this subsection with a lemma which will be used later on:

\lemmax{
\label{conginZ}
Let $\tau,\mu \in J$, $h \in \PL_+(J)$. Then:
\begin{itemize}
\item[(i)]The limit $\varphi_\pm=\underset{n \to \infty}{\lim}
h^{\pm n}(\tau)$ exists and $h(\varphi_\pm)=\varphi_\pm$,
\item[(ii)]We can determine whether there is or
not an $n \in \mathbb{Z}$, such that $h^n(\tau) = \mu$.
\end{itemize}
}
\begin{proof}
The two sequences
$\{h^{\pm n}(\tau)\}_{n \in \mathbb{N}}$
are strictly monotone, and they have a limit
$\underset{n \to \infty}{\lim} h^{\pm n}(\tau)= \varphi_\pm \in J$.
Thus, by continuity of $h$
\[
\varphi_\pm = \underset{n \to \infty}{\lim}h^{n+1}(\tau)=
\underset{n \to \infty}{\lim}h(h^n(\tau))= h(\varphi_\pm).
\]

Thus we have that
$\{h^{n}(\tau)\}_{n \in \mathbb{Z}} \subseteq (\varphi_-,\varphi_+)$
and we have that
$\varphi_{\pm}$ are the closest intersections of the graph of $h$
with the diagonal on the point $\tau$.
It is possible to compute $\varphi_+,\varphi_-$ directly,
since the graph of $h$ is piecewise linear.
As a first check, we must see if $\mu$ is between the points
$\varphi_-$ and $\varphi_+$.
Then since the two sequences $\{h^{\pm n}(\tau)\}_{n \in \mathbb{N}}$
are monotone, then after a finite number of steps we find
$n_1,n_2 \in \mathbb{Z}$ such that $h^{-n_1}(\tau) < \mu
< h^{n_2}(\tau)$ and so this means that either there is an integer
$-n_1 \le n \le n_2$ with $h^n(\tau)=\mu$ or
not, but this is a finite check.
\end{proof}

\subsection{The Stair Algorithm for $\PL^0_{S,G}(J)$}
\label{ssec:stairfull}

Section~\ref{ssec:diagonls} will prove that we can reduce our study to
$y$ and $z$ such that $\Fix(y) = \Fix(z)$.
Recall that an intersection point $\alpha $
of the graph of $z$ with the diagonal
may not be a point in $S$ (for instance,
a dyadic rational for the case of $\PL_2(I)$;
see again figure~\ref{fig:rational-intersection}).
If this is the case then $\alpha$ cannot be a breakpoint for $y,z$
and more importantly for $g$.
Recall that, by Definition~\ref{thm:define-PLSG-no-intersection},
a function $z$ is in $\PL_{S,G}^0(J)$ if $\Fix(z)$
does not contain any point of $S$, except for the endpoints of $J$.
\propositionx{
\label{thm:stairfull}
Let $y,z \in \PL^0_{S,G}(J)$ and and $q$ be a fixed element in $G$.
Suppose that $\Fix(y)=  \Fix(z)$.
We can decide whether or not there is a
$g \in \PL_{S,G}(J)$ with initial slope $g'(\eta^+)=q$
such that $y$ is conjugate to $z$ through $g$.
If $g$ exists it is unique.
Moreover there is an algorithm for constructing this conjugator.
}

\begin{proof}
This proof will be essentially the same as the previous stair algorithm
with a few more remarks.
We assume therefore that such a conjugator exists and build it. Let
$\Fix(y)= \Fix(z)=
\{\eta=\alpha_0  < \alpha_1 < \ldots < \alpha_s < \alpha_{s+1}=\zeta\}$.
We restrict our attention to $\PL_{S,G}([\alpha_i,\alpha_{i+1}])$
(as defined in Remark~\ref{thm:define-PLSG}), for each $i=0,\ldots,s$.
If $y$ and $z$ are conjugate on $[\alpha_i,\alpha_{i+1}]$
then we can speak of linearity
boxes: let $\Gamma_i:=[\alpha_i,\gamma_i] \times [\alpha_i,\gamma_i]$
be the initial linearity box and $\Delta_i:=[\delta_i,\alpha_{i+1}] \times
[\delta_i,\alpha_{i+1}]$ the final one for $\PL_{S,G}([\alpha_i,\alpha_{i+1}])$.
Now what is left to do is to repeat the procedure of the stair algorithm
for elements in
$\PL_{S,G}^<(U)$,
for some interval $U$. We build a conjugator $g$ on $[\alpha_0,\alpha_1]$
by means of the stair  algorithm.
We observe that $\alpha_1$ is not a breakpoint,
hence $g'(\alpha_1^+)=g'(\alpha_1^-)$.
Thus we are given an initial
slope for $g$ in $[\alpha_1,\alpha_2]$,
then we can repeat the same procedure and repeat the stair algorithm
on $[\alpha_1,\alpha_2]$. We keep repeating the same procedure until we
reach $\alpha_{s+1}=\zeta$. Then we check whether the $g$ we have found conjugates
$y$ to $z$.
Finally, we observe that in each square
$[\alpha_i,\alpha_{i+1}]\times[\alpha_i,\alpha_{i+1}]$
the determined function is unique,
since we can apply Lemma~\ref{cong:unique} on it.
\end{proof}

An immediate consequence of the previous result is the following Lemma:

\lemmax{
\label{thm:only-centralizer-identity}
Suppose $z \in \PL^0_{S,G}(J)$ and $g \in \PL_{S,G}(J)$ are such that
$g'(\eta^+)=1$ and $(g^{-1}zg)(t)=z(t)$, for all $t\in J$.
Then $g(t)=t$, for all $t \in J$.
}

\remarkx{
\label{backwardstairalgorithm}
It is possible to run
a backwards version of the stair algorithm also for $\PL_{S,G}^0(J)$.
Moreover, in this case it also possible to run a midpoint version of it:
if we are given a point $\lambda$ in the interior of $J$ fixed by $y$ and $z$
and $q \in G$ , then, by running
the stair algorithm at the left and the right of $\lambda$
we determine whether there is or not
a conjugator $g$ such that $g'(\lambda)=q$. }


\notationx{We recall that, given $y \in \PL_{S,G}(J)$,
we denote the centralizer of $y$ in $\PL_{S,G}(J)$ by
\[
C_{\PL_{S,G}(J)}(y)=\{g \in \PL_{S,G}(J) \mid y^g=y \}.
\]}

From Lemma~\ref{thm:only-centralizer-identity} and
Remark~\ref{backwardstairalgorithm} we have:

\corollaryx{
\label{con:unique}
Let $y,z \in \PL^0_{S,G}(J)$ such that $\Fix(y)=\Fix(z)$
and let
$$
C_{\PL_{S,G}(J)}(y,z)=\{g \in \PL_{S,G}(J) \, | \, y^g=z\}
$$
be the set of all conjugators. For any $\tau \in \Fix(y)$
define the map
$$
\begin{array}{lrcl}
\varphi_{y,z,\tau}: & C_{\PL_{S,G}(J)}(y,z) & \longrightarrow & \mathbb{R_+}  \\
               &  g               & \longmapsto     & g'(\tau),
\end{array}
$$
where if $\tau$ is an endpoint of $J$ we take only a one-sided derivative. Then
\begin{itemize}
\item[(i)]$\varphi_{y,z,\tau}$ is an injective map.
\item[(i)]If $\varphi_{y,z,\tau}$ admits a section, i.e.,
if there is a partially defined map $\mathbb{R_+} \to C_{\PL_{S,G}(J)}(y,z)$,
$\mu \to g_\mu$ such that $\varphi_{y,z,\tau}(g_\mu)=\mu$
then $g_\mu$ is unique and can be constructed.
\end{itemize}
}

\propositionx{
\label{thm:psi-surjective2}
Let $y,z \in \PL_{S,G}^0(J)$ such that $\Fix(y)=\Fix(z)$
and let $\lambda$ be in the interior of $J$ such that
$y(\lambda) \ne \lambda$. Define
$$
\begin{array}{lrcl}
\psi_{y,z,\lambda}:      & C_{\PL_{S,G}(J)}(y,z)   & \longrightarrow & J \\
                         & g                  & \longmapsto     & g(\lambda).
\end{array}
$$
Suppose $y^{n}(\lambda)\underset{n \to \infty}{\longrightarrow} \tau$. Then
\begin{itemize}
\item[(i)]There is a map $\rho_\lambda:J \to \mathbb{R}_+$
such that the following diagram commutes:
\begin{displaymath}
\begin{diagram}[grid=small]
                  &                                    & J\\
                  & \ruTo^{\psi_{y,z,\lambda}}         & \dTo_{\rho_\lambda} \\
C_{\PL_{S,G}(J)}(y,z) & \rTo_{\varphi_{y,z,\tau}}       & \mathbb{R}_+ 
\end{diagram}
\end{displaymath}
\item[(ii)]$\psi_{y,z,\lambda}$ is injective.
\item[(iii)]If $\psi_{y,z,\lambda}$ admits a section, i.e.,
if there is a partially defined map
$J \to C_{\PL_{S,G}(I)}(y,z)$,
$\mu \to g_\mu$ such that $\psi_{y,z,\lambda}(g_\mu)=\mu$
then $g_\mu$ is unique and can be constructed.
\end{itemize}
}

\begin{proof}
Let
$\Fix(y)=\Fix(z)=
\{\eta=\mu_0 < \mu_1 < \ldots < \mu_k < \mu_{k+1} = \zeta\}$
and suppose $\mu_i < \lambda < \mu_{i+1}$ for some $i$.
We define the partial map $\rho_\lambda:J \to \mathbb{R}_+$ as
\[
\rho_\lambda(\mu)=\begin{cases}
\lim_{n \to \infty} \frac{y^n(\mu)-\tau}{z^n(\lambda)-\tau} &
                            \mu \in [\mu_i,\mu_{i+1}] \\
1 & \mbox{otherwise}
\end{cases}
\]
Since $\Fix(y)=\Fix(z)$,
$z^{n}(\lambda)\underset{n \to \infty}{\longrightarrow} \tau$ and
$\tau$ is fixed by $g$. Thus if $\mu=g(\lambda)$,
then $y^n(\mu)=g(z^n(\lambda))\underset{n \to \infty}{\longrightarrow} \tau$.
With this definition, the proof follows closely that of
Lemma~\ref{cong:unique}(ii), Proposition~\ref{thm:psi-surjective}
and by applying Corollary~\ref{con:unique} and
Remark~\ref{backwardstairalgorithm}.
\end{proof}

\begin{figure}
 \begin{center}
  \includegraphics[height=5cm]{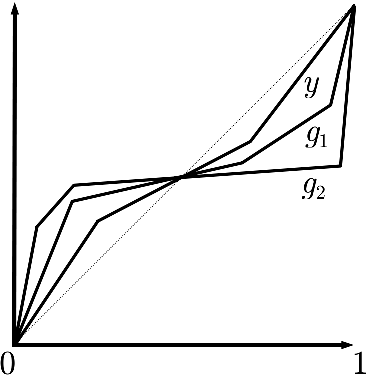}
 \end{center}
 \caption{Two centralizers $g_1,g_2$ of a function $y \in \PL^0_{S,G}(J)$.}
 \label{fig:intersect-fixed-points}
\end{figure}

Geometrically this says that if $y \in \PL^0_{S,G}$
then the graphs of the centralizers of $y$ inside $\PL^0_{S,G}$
intersect only at the fixed points of $y$ 
(see figure~\ref{fig:intersect-fixed-points}), which justifies the 
terminology ``almost one bump'' functions.

\section{Centralizers in subgroups of $\PL_+(I)$}
\label{sec:cenPL}

In this section we use the Stair algorithm to derive
several results about centralizers of elements in $\PL_+(J)$ and
$\PL_{S,G}(J)$. Although most of these results are already known,
our approach is new.
The main result of Subsection~\ref{ssec:centralizers-brin}
was first obtained by Brin and Squier~\cite{brin2}.
We will provide a new proof which generalizes to the case
$\PL_{S,G}(J)$.
The tools of Subsection~\ref{ssec:centralizers-brin} and the
results and proofs in the remaining subsections are new and constructive
(except for the results on the special case of Thompson's group $F$),
giving a procedure to solve the simultaneous conjugacy problem. We start by
giving an easy application of the stair algorithm before getting into the conjugacy problem.

\subsection{Centralizers of elements in $\PL^0_+(I)$ and $\PL_{S,G}^0(I)$\label{ssec:centralizers-brin}}

The Stair Algorithm from Section~\ref{sec:main-machine}
does not tell us anything about the image of
the homomorphism $\varphi_z : C_{\PL_+(J)}(z) \to \mathbb{R}_+$.
In this section we will show that if $z$ is
in $\PL_+^<(J)$ then the image of  $\varphi_z$
is a discrete subgroup of $\mathbb{R}^+$, thus
the centralizer of $z$ is an infinite cyclic group.
Let $A_z =\varphi_z \left(C_{\PL_+(J)}(z)\right) \subset \mathbb{R}_+$
be the set of all possible
initial slopes of centralizers.
The set $A_z$ is infinite, since $\langle z \rangle \subseteq C_{\PL_+(J)}(z)$.
Using the injectivity of $\varphi_z$,
we can define $\psi_z$ inverse of the function $\varphi_z$ on $A_z$
\[
\begin{array}{rcl}
\psi_z: A_z &\to& C_{\PL_+(J)}(z) \\
\alpha &\mapsto& g_\alpha,
\end{array}
\]
which is clearly a group isomorphism.
The previous section provides an algorithm to determine if
$c\in \mathbb{R}$ is an element in $A_z$ and
constructs the piecewise linear function $\psi_z(c)$ if it is defined.
which sends an initial slope $\alpha$ to its
associated conjugating function $g_\alpha$.

The main result of this section is the following:

\theoremx{
\label{thm:new-centralizer}
Let $J \subseteq [0,1]$ be a closed interval  and let
$id \ne z \in \PL_+^<(J)$. Then $C_{\PL_+(J)}(z)$ is isomorphic with $\mathbb{Z}$,
moreover there is an algorithm which constructs a generator $w$
of this group and $w$ is a root of $z$.
}

We remark that Theorem~\ref{thm:new-centralizer}
has originally been proved by Brin and Squier (Theorem 5.5 in~\cite{brin2}).
The connection between our proof and the one of Brin and Squier
has been described in a paper by the second author \cite{matucci3}.
We also observe that Altinel and Muranov gave another proof of this result using different methods
(Lemma 4.2 in~\cite{altimuran}).
The tools that we will use in our version of the proof
that we are about to give are relevant for
Lemma~\ref{thm:possible-interval-initial-slopes-PLSG},
which is central in our construction
of candidate conjugators.

\begin{proof}[Proof of Theorem~\ref{thm:new-centralizer}:]
By the discussion above we have that
the group 
$$
A_z=\{g'(\eta^+) \mid g \in C_{\PL_+(J)}(z)\}
$$
 is isomorphic with
$C_{\PL_+(J)}(z)$. We start by assuming that
$z \in \PL_+^<(J)$ and we want to prove that $A_z$ is discrete,
since any discrete subgroup of $\mathbb{R}_+$
is isomorphic to $\mathbb{Z}$. The argument below not
only proves that $A_z$ is discrete but also provides an
algorithm to find a generator of this group.

The proof relies on the following key lemmas:

\lemmax{
\label{non-lenarity}
Let $r$ be a positive integer such that $z^r({\beta}) < \alpha$, where
$[\eta,\alpha]^2$ and $[\beta,\zeta]^2$
are initial and final linearity boxes for the element $z$.
Either $z^{r}$ is not linear on the interval $[\beta, z^{-r}(\alpha)]$ or
$z^{2r}$ is not linear on the interval $[\beta, z^{-2r}(\alpha)]$.
}
\begin{proof}
Assume that both $z^r$ and $z^{2r}$ are both linear on these intervals
and denote their slopes by $s_1$ and $s_2$ respectively.
Using the linearity boxes for $z$ it can be seen that $z^r$
is linear on $[\eta,\alpha]$
with slope $a^r$, where $a=z'(\eta^+)$ and  $z^r$
is linear on $[z^{-r}(\beta),\zeta]$ with slope $b^r$, where $b=z'(\zeta^-)$.
Since $z^{2r} = z^r \circ z^r$ we get that $z^{2r}$
is linear on $[\beta, z^{-r}(\alpha)]$ with slope
$a^rs_1$ and also it is linear on  $[z^{-r}(\beta), z^{-2r}(\alpha)]$
with slope $b^rs_1$. Thus we have
$$
a^rs_1= s_2 = b^rs_1,
$$
however this is a contradiction because $a < 1 < b$ and $s_1 \not=0$.
\end{proof}

\lemmax{
\label{finitely_many_slopes}
Let $s$ be a positive integer such that $z^s({\beta}) < \alpha$ and
$z^{s}$ is not linear on the interval $[\beta, z^{-s}(\alpha)]$.
Then there exists $\varepsilon >0$ such that there are only finitely many
$g \in C_{\PL_+(J)}(z)$ such that
$1 \geq g'(\eta^+) \geq 1-\varepsilon$ and $1 \leq g'(\zeta^-) \leq 1+\varepsilon$
and they can be constructed.
}
\begin{proof}
Since $z^s$ is not linear on $[\beta, z^{-s}(\alpha)]$
there exists  $\varepsilon >0$ such that $z^s$ has breakpoints on
$I_\varepsilon=[(\beta+\varepsilon \zeta)/(1+\varepsilon), z^{-s}(\varepsilon\eta + (1-\varepsilon)\alpha)]$.
Let $\{\mu_1< \ldots < \mu_k \}$ be the breakpoints of $z^s$ on this interval.

For any $g \in C_{\PL_+(J)}(z)$ with
$1 \geq g'(\eta^+) \geq 1-\varepsilon$ and $1 \leq g'(\zeta^-) \leq 1+\varepsilon$
the linearity boxes give us that $g$ is linear on
$[\eta, \alpha]$ and $[(\beta+\varepsilon \zeta)/(1+\varepsilon), \zeta]$ and,
if $\varepsilon>0$ has been chosen small enough, the sets $I_\varepsilon$ and
$g^{-1}(I_\varepsilon)$ are not disjoint.
By construction, the breakpoints of $g\circ z^s$ on $I_\varepsilon$ are
$\{\mu_1< \ldots < \mu_k \}$
and the breakpoints of $z^s\circ g$ on $g^{-1}(I_\varepsilon)$ are
$\{g^{-1}(\mu_1)< \ldots < g^{-1}(\mu_k) \}$.
However for all but finitely many choices for
$g'(\zeta^-)$
the sets $\{\mu_1< \ldots < \mu_k \}$ and
$\{g^{-1}(\mu_1)< \ldots < g^{-1}(\mu_k) \}$ are disjoint.
Therefore $g\circ z^s \not= z^s\circ g$, which contradicts the assumption that
$g \in C_{\PL_+(J)}(z)$. Let $V \subseteq [1-\varepsilon,1]$
be the finite set of admissible final slopes $g'(\zeta^-)$ found before.
We run the backwards stair algorithm on each
slope in $V$ and determine which element centralizes $z$.
\end{proof}

Lemmas~\ref{non-lenarity} and~\ref{finitely_many_slopes}
immediately gives that $A_z$ is discrete,
which completes the proof of the first part of Theorem~\ref{thm:new-centralizer}.
To construct a generator $v$ for $C_{\PL_+(J)}(z)$, we observe that $v^k=z$
for some integer $k$, hence $v$ is a root of $z$ and so
$v'(\zeta^-) \in [z'(\zeta^-),1]$. Let $V \subseteq [1-\varepsilon,1]$
be the set of admissible final slopes for $g \in C_{\PL_+(J)}(z)$
given by Lemma~\ref{finitely_many_slopes}.
We run the backwards stair algorithm on the finite set
\[
\Big( [z'(\zeta^-),1-\varepsilon] \cup V\Big) \cap
\{\sqrt[m]{z'(\zeta^-)}\}_{m \in \mathbb{Z}}
\]
of admissible final slopes and
pick the centralizing element $w$ with initial slope closest to $1$.
By injectivity of the map $\varphi_{z}$ of Lemma~\ref{sec:unque},
the map $w$ is a generator for $C_{\PL_+(J)}(z)$.
\end{proof}

We finish with a generalization of Lemma~\ref{finitely_many_slopes}: the
following result is central in solving the
simultaneous conjugacy problem in $\PL_+(I)$ (together with
the stair algorithm (Corollary~\ref{thm:explicit-conjugator})
and Lemma~\ref{thm:initial-box}).
It provides an algorithm for restricting the initial slopes
of the conjugators. Not only this allows us to effectively
solve the conjugacy problem in $\PL_+(I)$ but also allows
us the extend this solution to the groups $\PL_{S,G}(J)$,
provided that the additive group $S$ satisfy some mild computational requirements.

\lemmax{
\label{thm:possible-interval-initial-slopes-PLSG}
Let $J = [\eta, \zeta]$ be a closed interval with endpoints in $S$
and let $c > 1$, then the set
\[
N = \{g \mid g \in C_{\PL_+(J)}(y,z), \, g'(\eta) \in[c^{-1},c],
g'(\zeta) \in[c^{-1},c] \}
\]
is finite and  can be constructed.}
\begin{proof}
Let $\alpha' = \eta + c^{-1}(\alpha -\eta)$ and
$\beta' = \zeta - c^{-1} (\zeta -\beta)$. Using Lemma~\ref{thm:initial-box}
we can see that $g \in N$
is linear on the intervals $[\eta,\alpha']$ and $[\beta',\zeta]$.
By lemma~\ref{non-lenarity} there exists a sufficiently big integer $s$ such that
$z^{-s}(\alpha') \geq \beta'$ and $z^{s}$ is not linear on the interval
$[\beta', z^{-s}(\alpha')]$. Let $\mu_i$ denote the set of breakpoints of $z^{s}$ on this interval.
The function $g z^{s}$ has $\mu_i$ as breakpoints  since $g$ is linear in the first linearity box.
If $g$ is a conjugator we have that $g z^{s}=y^{s}g$ therefore $\mu_i$ are breakpoints of $y^{s}g$,
which means that $g(\mu_i)$ are breakpoints of $y^{s}$. This condition leaves finitely many possibilities
for the final slope of $g$, which shows that the set $\{g'(\eta) | g \in N \}$ is finite.
For each of the slopes in $\{g'(\eta) | g \in N \}$ we can construct a candidate conjugator and
test it.
\end{proof}


\theoremx{
\label{thm:new-centralizer2}
Let $J \subseteq [0,1]$ be a closed interval with endpoints in $S$ and let
$id \ne z \in \PL_{S,G}^0(J)$. Then $C_{\PL_{S,G}(J)}(z)$ is isomorphic with $\mathbb{Z}$,
moreover there is an algorithm which constructs a generator $w$ of this group and $w$ is a root of $z$.
}
\begin{proof}
Let $\partial \Fix(z)=\{\eta < \gamma< \ldots <\zeta\}$ and consider the injective homomorphism $\lambda$
defined by sending each element of $C_{\PL_{S,G}(J)}(z)$ to its restriction to the interval $[\eta,\gamma]$.
By construction, the image of $\lambda$ is contained in $C_{\PL_+(J)}(z) \cong \mathbb{Z} = \langle w\rangle$, hence
$C_{\PL_{S,G}(J)}(z)$ is also infinite cyclic and contains $z$.
By~Lemma~\ref{finitely_many_slopes} there are only finitely many
admissible initial slopes to be tested, so to find a generator we follow the same procedure used in
Theorem~\ref{thm:new-centralizer}.
\end{proof}

\subsection{Centralizers of elements in $\PL_+(J)$ and $\PL_{S,G}(I)$}

The results about centralizers of elements in  $\PL_+^0(J)$ and $\PL_{S,G}^0(I)$
can be extended to arbitrary elements by observing that any centralizer of $y$
need to fix all points in $\partial_S \Fix(y)$.
\theoremx{
\label{thm:centralizers-PL+G}
Let $J = [\eta,\zeta] \subseteq [0,1]$ be a closed interval
and $z \in \PL_+(J)$, then:
\begin{itemize}
\item[(i)] $C_{\PL_+(I)}(z)$ is isomorphic with a direct product of copies of 
$\mathbb{Z}$
and $\PL_+(J_i)$ for some suitable intervals $J_i \subseteq I$.

\item[(ii)]For every positive integer $n$ we can decide whether or not $\sqrt[n]{z}$ exists.
The map $z$ has only a finite number of
roots and every root is constructible, i.e., there is an algorithm to compute it.
\end{itemize}
}
\begin{proof}
(i) Consider the conjugacy problem with $y=z$ and
let 
$$
\partial \Fix(z)=
\{\eta=\alpha_0  < \alpha_1 < \ldots < \alpha_s < \alpha_{s+1}=\zeta\}.
$$
Any centralizer $g$ of $z$ must fix the set $\partial \Fix(z)$
and thus each of the $\alpha_i$'s.
Therefore we compute the centralizer of the restrictions $z_i$ of $z$
in each of the subgroups $\PL_+(J_i)$, where $J_i=[\alpha_i,\alpha_{i+1}]$ and so we
can assume that $z_i \in \PL_+^<(J_i)$ or $z_i \in \PL_+^<(J_i)$ or $z_i=id$.
If $z_i=id$, then it is immediate that $C_{\PL_+(J_i)}(z_i)=\PL_+(J_i)$.
Suppose $z \ne id$ on $[0,1]$, then, by Theorem~\ref{thm:new-centralizer}, we have
$C_{\PL_+(J_i)}(z_i) \cong  \mathbb{Z}$.

(ii) Again we suppose that
$\partial \Fix(z) = \{ 0 = \alpha_0  < \alpha_1 < \ldots < \alpha_r < \alpha_{s+1} = 1\}$
and we restrict to an interval $[\alpha_i,\alpha_{i+1}]$. Let $m = z'(0^+)$. We want to determine whether or not there is
a $g \in \PL_+([\alpha_i,\alpha_{i+1}])$ such that $g^{-1}zg=z$ and such that $g'(0^+)=\sqrt[n]{m}$.
Suppose that there is such a $g$, then $g^{-n} z g^n = z$ and $(g^n)'(0^+)=m$.
By injectivity of the map $\varphi_{z,z,\alpha_i}$ (Corollary~\ref{con:unique}), we have that $g^n=z$. Conversely,
if we have $h$ such that $h^n=z$, then $h'(0^+)=\sqrt[n]{m}$. But $h^{-1}zh= h^{-1} h^n h = h^n =z$.
Thus an element $h$ is a $n$-th root of $z$ if and only if it
is the solution the ``differential type'' equation with a given initial condition
$$
\begin{cases}
h^{-1}zh=z \\
h'(0^+)=\sqrt[n]{m}.
\end{cases}
$$

So we can decide whether or not there is a $n$-th root, by solving the equivalent conjugacy problem
with a given initial slope.
Moreover, if the $n$-th root of $g$ exists,  it is computable by
Proposition~\ref{thm:stairfull} and unique by Corollary~\ref{con:unique}.
Moreover, only finitely many roots are possible: the sequence $\sqrt[n]{z'(\eta^+)}$ converges to $1$
but Lemma~\ref{finitely_many_slopes} implies that only finitely many elements of this sequence can be candidate slopes
for a root.
\end{proof}

\propositionx{
\label{thm:break-centralizer}
Let $x \in \PL_+(J)$ and $\alpha$ be a point  in $J$. If $g \in C_{\PL_+(I)}(x)$ and
$g(\alpha)=\alpha$ then the functions
$$
g_{<,\alpha} = \left\{ \begin{array}{ll}
t & \mbox{ if } t \leq \alpha\\
g(t) & \mbox{ if } t \geq \alpha
\end{array} \right.
\quad
\quad
g_{>,\alpha} = \left\{ \begin{array}{ll}
g(t) & \mbox{ if }\,\, t \leq \alpha\\
t & \mbox{ if }\,\, t \geq \alpha
\end{array} \right.
$$
are also in the centralizer $C_{\PL_+(J)}(x)$ and $g$ is equal to the product of
$g_{<,\alpha}$ and $g_{>,\alpha}$.
}
\begin{proof}
If $x(\alpha)=\alpha$, this follows from Theorem~\ref{thm:centralizers-PL+G}.
Assume now that $x(\alpha) \ne \alpha$ and let $[c,d]$
be the largest interval containing $\alpha$ on which $x$ is a one-bump function.
Since $g$ centralizes $x$, the points $c$ and $d$ are fixed by both $x$ and $g$ and,
in particular, the proposition follows for the maps $g_{<,c}$ and $g_{>,c}$.
The conclusion will then follow if we can prove that $g_{<,\alpha}=g_{<,c}$ and $g_{>,\alpha}=g_{>,c}$.
The restriction $g|_{[c,d]}$ centralizes $x|_{[c,d]}$ and so, by Theorem
\ref{thm:new-centralizer}, we have that $g|_{[c,d]}=(\sqrt[m]{x})^k$ for suitable integers $m,k$. Since $x(\alpha) \ne \alpha$,
then $k=0$ and $g|_{[c,d]}=id$. It is now straightforward to verify
that $g_{<,\alpha}=g_{<,c}$ and $g_{>,\alpha}=g_{>,c}$.
\end{proof}

We will see that solving the simultaneous conjugacy problem is equivalent to be
detecting whether or not a given candidate function lies in the
intersection of finitely many centralizers.
The next results shows that the intersection of centralizers has structure
similar to a single centralizer, which allows us to modify the solution
of the conjugacy problem in $\PL_+(J)$ and $\PL_{S,G}(J)$ and verify
is it is possible to find conjugator in the intersection of several centralizers.

\propositionx{
\label{intersectionofcetralizers-PL_+}
Let $x_1, \ldots, x_k \in \PL_+(J)$ and define
$C:=C_{\PL_+(J)}(x_1) \cap \ldots \cap C_{\PL_+(J)}(x_k)$.
If the interval $J$ is divided by the points in the union
$\partial \Fix(x_1) \cup \dots \cup \partial \Fix(x_k)$
into the intervals $J_i$ then
$$
C = C_{J_1} \cdot C_{J_2} \cdot \ldots \cdot C_{J_r},
$$
where $C_{J_i}:= \{f \in C \mid f(t)=t, \forall t \not \in J_i\} = C \cap \PL_+(J_i)$.
Moreover, each $C_{J_i}$ is isomorphic to either $\mathbb{Z}$, or $\PL_+(J_i)$ or
is the trivial group.}

\begin{proof}
The set $\partial \Fix(x_i)$ is fixed by all elements in $C_{\PL_+(J)}(x_i)$,
therefore all elements in $C$ fix the end points of the intervals $J_i$,
since for $\alpha \in \cup \partial \Fix(x_i)$ and any $g\in C$ the function
$g_{<,\alpha}$ and $g_{>,\alpha}$ are in $C$ by
Proposition~\ref{thm:break-centralizer}.
Any element $z \in C$ can be written as the product
$z_1 \ldots z_r$, where $z_i \in \PL_+(J)$ is trivial outside of $J_i$
and $z_i|_{J_i} \in C_{\PL_+(J_i)}(x_n|_{J_i})$,
for all $n=1,\ldots, r$. Hence $z_i \in C_{J_i}$.
\end{proof} 

\corollaryx{
\label{thm:intersection-two-PL_+}
The intersection of any number $k\ge2$ centralizers of elements $x_1,\ldots,x_k$
in $\PL_+(J)$ is equal to the intersection of centralizers of two elements
$w_1,w_2 \in \PL_+(J)$
which are not necessarily part of the initial set $\{x_1,\ldots,x_k\}$.
}
\begin{proof}
Let $C=C_{\PL_+(I)}(x_1) \cap \ldots \cap C_{\PL_+(I)}(x_k)$
be the intersection of $k\ge 2$ centralizers of elements of $\PL_+(J)$.
By the previous Proposition we have $I=J_1 \cup \ldots \cup J_r$
and $C=C_{J_1} \cdot \ldots \cdot C_{J_r}$.
We want to define $w_1,w_2 \in \PL_+(I)$ such that
$C=C_{\PL_+(I)}(w_1) \cap C_{\PL_+(I)}(w_2)$.
We define them on each interval $J_i:=[\alpha_i,\alpha_{i+1}]$,
depending on $C_{J_i}$. \emph{Case 1:} If $C_{J_i}=id$,
then we define $w_1,w_2$ as any two elements in $\PL_+^<(J_i)$
so that one is not a power of the other.
\emph{Case 2:} If $C_{J_i}\cong \langle x \rangle$ for some
$id \ne x \in \PL_+(J_i)$,
then we define $w_1=w_2=x$.
\emph{Case 3:} If $C_{J_i}=\PL_+(J_i)$, then we define $w_1=w_2=id$.
\end{proof}


Using Theorem~\ref{thm:new-centralizer2} one can easily generalize the
results in the previous subsection
to the groups $\PL_{S,G}(J)$.

\theoremx{
\label{thm:centralizers-PLSG}
Let $J = [\eta,\zeta] \subseteq [0,1]$ be a closed interval with endpoints in $S$
and $z \in \PL_{S,G}(J)$, then:
\begin{itemize}
\item[(i)] $C_{\PL_{S,G}(I)}(z)$ is isomorphic with
a direct product of copies of the group $\mathbb{Z}$'s
and $\PL_{S,G}(J_i)$'s for some suitable intervals $J_i \subseteq I$.
\item[(ii)]For every positive integer $n$ we can decide whether or not
$\sqrt[n]{z}$ exists.
The map $z$ has only a finite number of
roots and every root is constructible, i.e., there is an algorithm to compute it.
\end{itemize}
}

\begin{proof}
(i) We consider the conjugacy problem with $y=z$ and let
$$
\partial_S \Fix(z)=
\{\eta=\alpha_0  < \alpha_1 < \ldots < \alpha_s < \alpha_{s+1}=\zeta\}.
$$
Any centralizer $g$ of $z$ must fix $\partial_S \Fix(z)$ pointwise.
We thus compute the centralizer of
the restrictions $z_i$ of $z$ in each of the
groups $\PL_{S,G}(J_i)$ where $J_i=[\alpha_i,\alpha_{i+1}]$
and assume that $z_i \in \PL_{S,G}^0(J_i)$ or $z=id$.
The rest of the proof follows as in Theorem~\ref{thm:centralizers-PL+G}(i)
by means of Theorem~\ref{thm:new-centralizer2}.

(ii) This follows as in Theorem~\ref{thm:centralizers-PL+G}(ii).
\end{proof}

Knowing the structure of a centralizer in $\PL_{S,G}(I)$ allow us the
extent the results about intersections of centralizers.

\propositionx{
\label{thm:intersectionofcentralizers}
Let $J = [\eta,\zeta] \subseteq [0,1]$ be a closed interval with endpoints in $S$,
let $z_1, \ldots, z_k \in \PL_{S,G}(J)$ and define the subgroup
$C:=C_{\PL_{S,G}(I)}(z_1) \cap \ldots \cap C_{\PL_{S,G}(I)}(z_k)$.
If the interval $J$ is divided by the points in the union
$\partial_S \Fix(z_1) \cup \dots \cup \partial_S \Fix(z_k)$
into the intervals $J_i$ then
$$
C = C_{J_1} \cdot C_{J_2} \cdot \ldots \cdot C_{J_r},
$$
where
$C_{J_i}:= \{f \in C \mid f(t)=t, \forall t \not \in J_i\} =
C \cap \PL_{S,G}(J_i)$.
Moreover, each
$C_{J_i}$ is isomorphic to either $\mathbb{Z}$, or $\PL_{S,G}(J_i)$ or
the trivial group.
}

\corollaryx{
\label{thm:intersection-two-PLSG}
The intersection of any number $k\ge2$ centralizers of elements $x_1,\ldots,x_k$
in $\PL_{S,G}(J)$ is equal to the intersection of centralizers of two elements
 $w_1,w_2 \in \PL_{S,G}(J)$
which are not necessarily part of the initial set $\{x_1,\ldots,x_k\}$.

}

\questionx{
Corollary~\ref{thm:intersection-two-PLSG} determines that any intersection
any number $k \ge 2$ centralizers elements
$x_1,\ldots,x_k$ in $\PL_{S,G}(J)$ can be expressed as the intersection
$C_{\PL_{S,G}(J)}(w_1) \cap C_{\PL_{S,G}(J)}(w_2)$ for two suitable elements
$w_1,w_2 \in \PL_{S,G}(J)$.
Is it possible to build the two elements $w_1,w_2$ inside the subgroup
$\langle x_1,\ldots,x_k\rangle$?
}


The groups $\PL_{S,G}(J_i)$ may not be isomorphic to each other
(see Remark~\ref{thm:multiple-orbits-of-points}).
However, in the special case of $S= \mathbb{Z}{\left[\frac{1}{2}\right]}$
it is true that $\PL_{S,G}(J_i) \cong F$, for
all $i$ (see Remark~\ref{thm:endpoints-in-S}).
This simplifies the statement of Theorem~\ref{thm:centralizers-PLSG} in
the case of Thompson's group $F$.
Also the proof can be simplified because one can use the
discreteness of the group $G$ instead of Lemmas~\ref{non-lenarity}
and~\ref{finitely_many_slopes} and Theorem~\ref{thm:new-centralizer2}.
As we have already mentioned this result is well known
and was first proved by Guba and Sapir~\cite{gusa1}
using different techniques.

\begin{figure}
 \begin{center}
  \includegraphics[height=5cm]{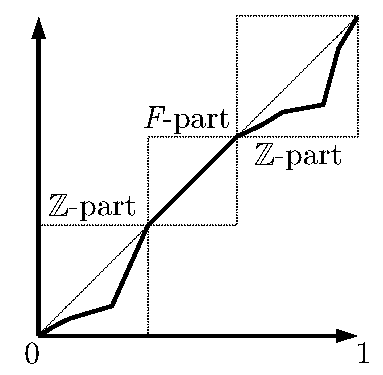}
 \end{center}
 \caption{The structure of centralizers in $F$}
 \label{fig:structure-centralizers-F}
\end{figure}

\theoremx{
\label{thm:centralizers-roots-thompson}
Let $z \in F \cong\PL_2(I)$. Then:
\begin{enumerate}
\item[(i)] Its centralizer is
$C_F(z) \cong F^m \times \mathbb{Z}^n$, for some positive integers $m,n$ such that
$0\le m \le n+1$ (see figure~\ref{fig:structure-centralizers-F}).
\item[(ii)] If $z \ne id$, the function $z$ has only a finite number of
roots and every root is constructible, i.e., there is an algorithm to compute it.
\end{enumerate}
}

\section{Moving fixed points}
\label{ssec:diagonls}

In this Section we describe Step 1 of the outline of
 Subsection~\ref{ssec:outline-strategy}. If two maps
$y,z$ are conjugate via $g$, then $g(\Fix(y))=\Fix(z)$.
Thus, moving fixed points is an intermediate step towards the conjugacy problem.
We begin our proofs for the easier case of $\PL_+(J)$
and then move on to study the case of the groups $\PL_{S,G}(J)$.

\subsection{Moving fixed points in $\PL_+(J)$}

This case is the easiest one --
essentially, in the case of $\PL_+(J)$, the only necessary thing to check is if
$\Fix(y)$ and $\Fix(z)$ have the same number and
``type'' of components and they have the ``same order''%
\footnote{This is exactly the invariant $\Sigma_2$ defined by Brin and Squier in \cite{brin2}.}.
We state without proof the following results:

\theoremx{
\label{movingpoints inPL_+}
Let $y_1 < y_2 < \dots < y_n$ and $z_1 < z_2 < \dots < z_n$ be points in the interval $J$.
Then there exists a $g \in \PL_+(J)$ such that $g(y_i)=z_i$ for all $i=1,\dots,n$.
}

\theoremx{
\label{congdiagonals-PL_+}
Let $y ,z \in \PL_+(J)$. There is a algorithm, which constructs
an element  $g \in \PL_+(J)$ such
that $g(\Fix(y)) = \Fix(g^{-1}yg)=\Fix(z)$,
or shows that such element does not exist.
}

\subsection{Moving fixed points in $\PL_{S,G}(J)$}

The main difference  between the groups $\PL_{S,G}(J)$ and $\PL_+(J)$ is
that (in general) $\PL_{S,G}(J)$ does not act transitively on the interior
points in the interval $J$. Our first step it to describe the orbits.
Let us define an equivalence relation $\sim_{S,G,J}$ in $J$. If
$x,y \in J$ we say that $x \sim y$
if and only if there exists $g\in \PL_{S,G}(J)$ such that $g(x)=y$.
Unless otherwise stated, we always assume that the endpoints of $J$ are in $S$.

\definitionx{Let $\mathcal{I}_{S,G}$ 
denote the submodule of the $\mathbb{Z}[G]$-module $S$
generated by $ (g-1)$ for $g \in G$.
We denote with $\pi_{S,G}:S \to S/\mathcal{I}_{S,G}$ the natural quotient map.
Unless otherwise stated, we will drop the subscript and
write $\mathcal{I}$ and $\pi$ instead of $\mathcal{I}_{S,G}$ and $\pi_{S,G}$.
\label{thm:ideal-def}}

We remark that the natural map $\pi$ is a homomorphism.
The next theorem plays central role in understanding the orbits of points in $J$
under the action of $\PL_{S,G}(J)$
by detecting when two points of $S$ are in the the same $\PL_{S,G}$-orbit.

\theoremx{
\label{moving_points_ideal}
Let $J$ be an interval with end points in $S$ and let $x,y \in S \cap J$.
Then $x\sim y$ if and only if $x-y \in \mathcal{I}$.
}

The proof follows from the next two results.

\lemmax{
\label{thm:nec-cond-trans}
Let $J \subseteq [0,1]$ be a closed interval with
at least one of the endpoints $\eta$ in $S$ and let $g \in \PL_{S,G}(J)$.
Then, for every $t \in J \cap S$, we have $\pi(g(t))=\pi(t)$.
}

\begin{proof}
We can assume that the $\eta$ is the left one
and we apply induction on the number of breakpoints preceding $t$. In case the endpoint in $S$
is the right one, we apply induction on the breakpoints following $t$. Let $\{\eta_1, \ldots, \eta_r\}$
be the set of all breakpoints of $g$ on the interval $[\eta,t)$.
Then $g(t)=c_r (t-\eta_r) +g(\eta_r)$ for some suitable $c_i \in G$.
By induction hypothesis, the number of breakpoints preceding $\eta_r$ is $r-1$ and
so we have $\pi(g(\eta_r))=\pi(\eta_r)$. Now we observe that
\begin{eqnarray*}
\pi(g(t))=\pi(c_r (t-\eta_r)+ g(\eta_r))= &\\
\pi(c_r -1)\pi(t-\eta_r) + \pi(1)\pi(t-\eta_r) + \pi(g(\eta_r))= &\\
\pi(t-\eta_r)+\pi(\eta_r)=\pi(t). &
\hfill \qedhere
\end{eqnarray*}
\end{proof}

\propositionx{
\label{thm:equiv-cond-trans}
Let $J \subseteq [0,1]$ be a closed interval with both endpoints in $S$
and let $u,v \in J \cap S$. Then $\pi(u)=\pi(v)$
if and only if there is a $g \in \PL_{S,G}(J)$ such that
$g(u)=v$.}

\begin{proof}
The sufficient condition is implied by Lemma~\ref{thm:nec-cond-trans}.
Suppose now that $J=[\eta,\zeta]$ and let $L=\zeta-\eta$.
We recenter the axis at $(\eta,\eta)$, so that interval $J$ is now $[0,L]$.
For $\alpha \in G, \beta \in J \cap S$ such that $\alpha \beta < L-\beta$
define (see figure~\ref{fig:transitive})
\[
g_{\alpha,\beta}(t):=\begin{cases}
\alpha t & t \in [0,\beta] \\
t- (1-\alpha)\beta & t \in [\beta, L-\alpha\beta] \\
\frac{1}{\alpha}(t-L) + L  & t \in [L-\alpha\beta,L]
\end{cases}.
\]

\begin{figure}[0.5\textwidth]
 \begin{center}
  \includegraphics[height=5cm]{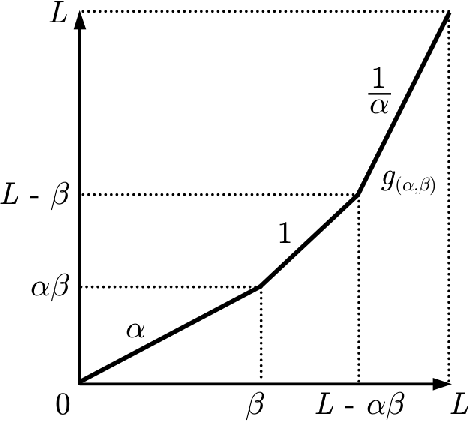}
 \end{center}
 \caption{The basic function to get transitivity.}
 \label{fig:transitive}
\end{figure}

Using the maps $g_{(\alpha,\beta)}$ or $g_{(\alpha,\beta)}^{-1}$ we can send any number $\beta \le t \le L-\alpha\beta$
to $t - (1-\alpha)\beta$ and any number $\alpha \beta \le t \le L -\beta$ to $t + (1-\alpha)\beta$.

Since $\pi(u)=\pi(v)$ then $v-u \in \mathcal{I}$ and so
\[
v - u =  (1-\alpha_1)\beta_1+ \ldots + (1-\alpha_k)\beta_k
\]
for some $\alpha_i \in G,\beta_i \in J \cap S$. Adding extra terms if necessary we can assume that
$$
u + (1-\alpha_1)\beta_1+ \ldots + (1-\alpha_i)\beta_i \in J
$$
for any $1\leq i \leq k$. Since $S$ is a dense subgroup of $\mathbb{R}$ then,
for each $\beta_i$, we can find numbers $\beta_{i,j} \in J\cap S$ small enough
such that 
\begin{itemize}
\item $L-\beta_{i,j}>\alpha_i\beta_{i,j}$ so that the map $g_{(\alpha_i,\pm \beta_{i,j})}$
can be defined, and
\item $\beta_i=\sum_{j} \beta_{i,j}$.
\end{itemize}
Finally we can see that the composition of the maps
$g_{(\alpha_i,\beta_{i,j})}^{\pm 1}$ sends $u$ to $v$, which finishes the proof.
\end{proof}

\corollaryx{
Any linear piece of the graph of a $g\in \PL_{S,G}(J)$ has equation of the form $x\to ax+b$
where $a\in G$ and $b\in \mathcal{I}$.
}

\corollaryx{
\label{two_intervals}
Let $J_1$ and $J_2$ be two intervals containing $x,y$, then $x \sim_{S,G,J_1} y$
if and only if  $x \sim_{S,G,J_2} y$.
}

\theoremx{
\label{thm:pre-extension}
Let $J$ be a closed interval with endpoints in $S$ and suppose we have
$u_1,v_1,\cdots,u_k,v_k \in J \cap S$ such that $u_1 < \cdots <u_k$, $v_1<\cdots<v_k$ and
$u_i \sim v_i$ for all $i=1,\ldots,k$.
Then there exists a $g \in \PL_{S,G}(J)$ such that $g(u_i)=v_i$ for all $i=1,\ldots,k$.
}
\begin{proof}
The proof is by induction: the base case $k=1$ is just the definition of the equivalence relation $\sim$.
Let $k>1$, by the induction assumption, there exist $ \widehat{g}\in \PL_{S,G}(J)$ such that
$\widehat{g}(u_i) = v_i$ for $i=1,\ldots,k-1$. Using that $\sim_{S,G,J}$ is an equivalence relation
we can obtain that $\widehat{g}(u_k)\sim_{S,G,J} u_k \sim_{S,G,J} v_k$. Let $J'$ denote the interval
$[v_{k-1}, \zeta]$ which contains the points $\widehat{g}(u_k)$ and $v_k$. By Corollary~\ref{two_intervals}
we have $\widehat{g}(u_k)\sim_{S,G,J'} v_k$, therefore there exists $\bar g \in \PL_{S,G}(J')$ such that
$\bar g (\widehat{g}(u_k)) = v_k$, thus the element $g = \bar g \circ \widehat{g}$ sends $u_i$ to $v_i$ for all $i$.
\end{proof}

\lemmax{
\label{extension-partial-maps}
Suppose $I_1,\ldots,I_k $ is a family of disjoint closed
intervals  $I_i=[a_i,b_i]$, with $b_i < a_{i+1}$
for all $i=1,\ldots,k$ and $a_i,b_i \in S$. Let
$J_1, \ldots, J_k \subseteq [0,1]$, with $J_i=[c_i,d_i]$,
be another family of intervals with the same property
such that $a_i \sim c_i$ and $b_i\sim d_i$.
Suppose that $g_i:I_i \to J_i$ is a piecewise-linear function with a finite number
of breakpoints, occurring at $S$
and such that all slopes are in $G$.
Then there exists a $\widetilde{g} \in \PL_{S,G}(I)$ such that $\widetilde{g}|_{I_i}=g_i$.
}
\begin{proof}
By Theorem~\ref{thm:pre-extension} there exists an $h \in \PL_{S,G}(J)$ with
$h(a_i)=c_i$ and $h(b_i)=d_i$. Define
\[
\widetilde{g}(t):=
\begin{cases}
h(t) & t \not \in I_1 \cup \ldots \cup I_k \\
g_i(t) & t \in I_i.
\end{cases}
\]
By construction, it is clear that $\widetilde{g} \in \PL_{S,G}(J)$ and $\widetilde{g}|_{I_i}=g_i$.
\end{proof}

\corollaryx{
Any part of the graph of $x\to ax+b$ where $a\in G$ and $b\in \mathcal{I}$,
inside the open square $J \times J$ can be extended to a graph
of an element in $\PL_{S,G}$.
}

Any isolated fixed point $\alpha$ of an element $g\in \PL_{S,G}(J)$ is of the form
$\alpha = s/(t-1)$ for some $s\in S$ and $t\in G \setminus \{1\}$. 
Let $Q_S$ denote the set of all points of the form $s/(t-1)$. 
The next step is to understand when two points in $Q_S$ are in one and the
same orbit under $\PL_{S,G}(J)$.

\theoremx{
\label{moving rational points}
Let $J=[\eta,\zeta]$ be a closed interval with endpoints in $S$ and
let $\alpha,\beta \in J \cap Q_S$.
The points $\alpha$ and $\beta$ are equivalent under $\sim_{S,G,J}$
if and only if we can find $s,s' \in S$ and $t\in G$ such that 
$\alpha=s/(t-1),\beta=s'/(t-1)$ and
\[
sG=s'G \pmod{(t-1) \mathcal{I}}
\]
where $(t-1) \mathcal{I}$ 
denotes the image of the submodule $\mathcal{I}$
under the multiplication by $t-1 \in \mathbb{Z}[G]$.
}

\begin{proof}
Suppose there is a map $g \in \PL_{S,G}(J)$ with $g(\alpha)=\beta$
and let $g(x)=cx+d$ in a small neighborhood $J_\alpha$ of $\alpha$.
We can choose representatives $s\in S$ and $t \in G$ such that
$\alpha=s/(t-1)$ and then, 
since $g \in \PL_{S,G}(J)$, we use Lemma~\ref{thm:nec-cond-trans} to get
\[
\pi(x)=\pi(g(x))=\pi(c-1)\pi(x) + \pi(x) + \pi(x)
\]
for all $x \in J_\alpha \cap S$ and therefore $\pi(d)=0$, which
implies $d \in \mathcal{I}$.
The equality $g(\alpha) = \beta$ implies that
$\beta = s'/(t-1)$ where $s' = cs + d(t-1)$, which implies that
$sG=s'G \pmod{(t-1) \mathcal{I}}$.

Conversely, suppose that we can write
$\alpha=s/(t-1),\beta=s'/(t-1)$, for some $s,s' \in S$ and $t\in G$ such that
 $sG=s'G \pmod{(t-1) \mathcal{I}}$.
The second condition implies that there exist $c_1,c_2 \in G, d_2 \in \mathcal{I}$ 
such that
\[
c_1 s = c_2 s' + (t-1) d_2
\]
and so if we set $c=c_2/c_1$ and $d=d_2/c_1$, we get $\alpha=c\beta +d$. 
Let $f(t)=ct+d$ be a line through the point $(\alpha,\beta)$ and let 
$[\gamma,\delta] \subseteq J$ be a small
interval such that $\gamma,\delta \in S$. 
Since $\pi(d)=0$ we have that $\pi(f(\gamma))=\pi(\gamma)$
and $\pi(f(\delta))=\pi(\delta)$ and so,
by Lemma~\ref{extension-partial-maps}
there is an $g \in \PL_{S,G}(J)$ with $g|_{[\gamma,\delta]}=f$.
By construction $g(\alpha)=\beta$ as required.
\end{proof}

Using the previous 2 results one can easily generalize
Theorem~\ref{congdiagonals-PL_+} to the groups $\PL_{S,G}(J)$.
Of course this is only possible if the group $S$ and the group
$G$ satisfy some mild computational requirements,
which are described in section~\ref{sec:comp_req}.

\corollaryx{Assume that $S$ and $G$ satisfy the computational requirements
from section~\ref{sec:comp_req}. Then for any $\alpha,\beta \in Q_S \cap J$
there is an algorithm which constructs to
$g \in \PL_{S,G}(J)$ such that $g(\alpha)=\beta$,
or shows that such element does not exist.
\label{thm:how-to-identify-rationsls}}

We state the same result for a finite number of points.
Its proof uses~Lemma~\ref{extension-partial-maps} on a number of
disjoint intervals, one around each point.

\begin{cor}
\label{thm:overlap-intersection-diagonal}
Assume that $S$ and $G$ satisfy the computational requirements
from section~\ref{sec:comp_req}.
Let $\eta<\alpha_1 < \ldots < \alpha_r<\zeta$ and $\eta<\beta_1 < \ldots < \beta_r<\zeta$ be two
partitions of $J$ with elements of the set $Q_S$.
Then there is an algorithm which constructs
$g \in \PL_{S,G}(J)$ with $g(\alpha_i)=\beta_i$,
or shows that such element does not exist.
\end{cor}

\theoremx{
\label{congdiagonals}
Assume that $S$ and $G$ satisfy the computational requirements
from section~\ref{sec:comp_req}.
Then given any $y,z \in \PL_{S,G}(I)$, there is an algorithm which constructs
$g \in \PL_{S,G}(I)$ such that
$g(\Fix(y)) = \Fix(g^{-1}yg)=\Fix(z)$,
or shows that such element does not exist.
}
\begin{proof}
First we check if $\# \partial \Fix(y) = \# \partial \Fix(z)$. Then we use the previous
Corollary to find a $g \in \PL_2(I)$, with $g(\partial \Fix(y))=\partial \Fix(z)$, if it exists.
To finish we check whether $\Fix(g^{-1}yg)$ contains the same intervals as $\Fix(z)$.
\end{proof}

\subsection{The case of Thompson's group}

Here are the analogues of previous results in the case of
Thompson's groups $F$.%
\footnote{The first two results are well known, see~\cite{cfp}.}

\lemmax{
\label{thm:standard-cannon-floyd-parry}
If $0=x_0< x_1 < x_2 < \ldots < x_n=1$
and $0=y_0 < y_1 < y_2 < \ldots < y_n=1$ are two partitions
of $[0,1]$ consisting of dyadic rational numbers, then we can build a $g \in \PL_2(I)$, such that $g(x_i)=y_i$.
}

An easy well known consequence is the following extension Lemma:

\lemmax{
\label{thm:extension-partial-maps}
Suppose $I_1,\ldots,I_k \subseteq[0,1]$ is a family of disjoint closed
intervals  $I_i=[a_i,b_i]$, with $b_i < a_{i+1}$
for all $i=1,\ldots,k$ and $a_i,b_i \in \mathbb{Z}[\frac{1}{2}]$. Let
$J_1, \ldots, J_k \subseteq [0,1]$, with $J_i=[c_i,d_i]$,
be another family of intervals with the same property.
Suppose that $g_i:I_i \to J_i$ is a piecewise-linear function with a finite number
of breakpoints, occurring at dyadic rational points,
and such that all slopes are integral powers of $2$.
Then there exists a $\widetilde{g} \in \PL_2(I)$ such that $\widetilde{g}|_{I_i}=g_i$.
}

\propositionx{
\label{rationals}
Let $\alpha=\frac{2^t m}{n}$ and $\beta=\frac{2^k p}{q}$ be rational numbers in $\mathbb{Q}\cap(0,1)$,
where $t,k \in \mathbb{Z}$, $m,n,p,q$ odd integers such that $(m,n)=(p,q)=1$.
Then there is a $g \in \PL_2(I)$ such that $g(\alpha)=\beta$ if
and only if $n=q$ and
\begin{equation}
\numberwithin{equation}{section}
\label{eq:prob}
p \equiv 2^Rm \pmod{n}
\end{equation}
for some $R \in \mathbb{Z}$. Equivalently there exist integers $t', k'$ such that
$2^{t'}\alpha - 2^{k'}\beta$ is an integer.
Moreover, there is an algorithm which constructs such element $g$
if the above condition is satisfied.
}

\examplex{
Let $\alpha=\frac{1}{17}$, $\beta=\frac{13}{17}$ and $\gamma=\frac{3}{17}$.
It is easy to see that we can find a $g \in \PL_2(I)$ with $g(\alpha)=\beta$, but there is no
$h \in \PL_2(I)$ with $h(\alpha)=\gamma$.}

We now state the analogue of Theorem~\ref{congdiagonals}
noticing that for Thompson's group the requirements
section~\ref{sec:comp_req} are satisfied.

\theoremx{
Given $y,z \in \PL_2(I)$, there is an algorithm which constructs
$g \in \PL_2(I)$ such that
$g(\Fix(y)) = \Fix(g^{-1}yg) = \Fix(z)$,
or shows that such element does not exist.
}

\section{The Conjugacy Problem and the Power Conjugacy Problem in $\PL_{+}(J)$ and $\PL_{S,G}(J)$}
\label{sec:applications}

The results of Section~\ref{ssec:diagonls}, together with the assumption that $S,G$
satisfy the computational requirements in section~\ref{sec:comp_req}, 
allow us to reduce to the problem to the case where
$\partial_S \Fix(y) = \partial_S \Fix(z)$.

\subsection{Characterizing Conjugacy in $\PL_+(J)$}

To study conjugacy between two elements $y$ and $z$ we can assume $\partial \Fix(y) =
\partial \Fix(z)= \{\alpha_1,\ldots,\alpha_n\}$ and we look for conjugators in $\PL_+([\alpha_i,\alpha_{i+1}])$
of the restrictions of $y$ and $z$ to $[\alpha_i,\alpha_{i+1}]$.
We reduce the study of the
conjugacy problem to smaller intervals. If $y=z=id$ on the interval $[\alpha_i,\alpha_{i+1}]$
there is nothing to prove, otherwise $y$ and $z$ are one-bump functions.
Given two elements $f,g \in \PL_+(J)$ we say that they are \emph{$y$-equivalent}
if $f=y^n g$, for some integer $n$.

\lemmax{\label{thm:conjugator-by-power} If $g$ is a conjugator for $y$ and $z$,
any $y$-equivalent map $y^n g$ is also a conjugator.}

\begin{proof} We observe that
\[
\left(y^{n} g\right)^{-1} y \left(y^{n} g\right)=g^{-1}yg=z.
\]
\end{proof}

\lemmax{\label{thm:bound-slopes-PL_+}
If $y,z \in \PL_+^<(J)$ are conjugate, there exists a $y$-equivalent conjugator $g \in \PL_+(J)$
such that $y(\lambda) < g(\lambda) < \lambda$, for any fixed $\lambda$ in the interior of $J$.}

\begin{proof}
Let $h \in \PL_+(J)$ be a conjugator for $y$ and $z$. Since $y \in \PL_+^<(J)$, there exists an integer $n$
such that $y^n h(\lambda) < y(\lambda) \le y^{n-1} h(\lambda)$.
By applying $y^{-1}$ on the inequality $y^n h(\lambda) < y(\lambda)$ we obtain
\[
y^n h(\lambda) < y(\lambda) \le y^{n-1} h(\lambda)<\lambda.
\]
We define $g=y^{n-1} h$ and we are done by Lemma~\ref{thm:conjugator-by-power}.
\end{proof}

\propositionx{
\label{thm:only-finitely-conjugators-tested}
To detect whether or not two elements $y,z \in \PL_+(J)$ are conjugate, only finitely many
functions need to be tested as possible candidate conjugators and they can be constructed.
Moreover we can enumerate all possible conjugators.}

\begin{proof}
By the discussion at the beginning of this subsection, we can assume that $y,z \in \PL_+^<(J)$.
Let $\lambda \in J$ be a fixed interior point of $J$ contained in the initial linearity box.
For any conjugator of $y$ and $z$, Lemma~\ref{thm:bound-slopes-PL_+} implies that
there is a $y$-equivalent conjugator $g \in \PL_+(J)$ such that $y(\lambda) < g(\lambda) <\lambda$.
Now, since the map $\rho_\lambda$ defined in Lemma~\ref{cong:unique} is increasing, it is immediate to see from
its definition that
\[
y'(\eta^+) = \rho_\lambda(y(\lambda)) \le g'(\eta^+) = \rho_\lambda(g(\lambda)) \le 1 =\rho_\lambda(\lambda) \le y'(\eta^+)^{-1}.
\]
Choosing another interior point $\mu$ in the final linearity box, we can use the analogous
version of $\rho_\mu$ at the final slope to obtain
$y'(\zeta^+)^{-1}\le g'(\zeta^+) \le y'(\zeta^+)$. Hence, the set of all conjugators $g$
such that $y(\lambda) < g(\lambda) <\lambda$ is contained in the set
\[
N:= \{h \mid h \in C_{\PL_+(J)}(y,z), \, h'(\eta) \in[y'(\eta^+),y'(\eta^+)^{-1}], h'(\zeta) \in[y'(\zeta^+)^{-1},y'(\zeta^+)] \},
\]
which by Lemma~\ref{thm:possible-interval-initial-slopes-PLSG} is finite and can be constructed.
If the set $N$ is non-empty then,
by the uniqueness of conjugators with a given initial slope (Lemma~\ref{cong:unique}) and by Lemma~\ref{thm:bound-slopes-PL_+},
the set of all conjugators for $y$ and $z$ is given by $\{y^r g \mid g \in N, r\in \mathbb{Z} \}$.
\end{proof}

\subsection{Conjugacy problem in $\PL_{S,G}(J)$}

We can now solve the conjugacy problem for elements in $\PL_{S,G}^0(J)$. We recall
that $\PL_{S,G}^0(J) \subseteq \PL_{S,G}(J)$ is the set of functions $f \in \PL_{S,G}(J)$ such that
the set $\Fix(f)$ does not contain elements of $S$ other than the endpoints of $J$.

\lemmax{
\label{thm:idea-conjugacy}
For any $y,z \in \PL_{S,G}^0(J)$ such that $y \ne z$ and $\Fix(y) = \Fix(z)$,
we can decide whether there is (or not) a $g \in \PL_{S,G}(J)$ with $y^g=z$.
Moreover we can construct and enumerate all possible conjugators.
}
\begin{proof}
In order to be conjugate, we must have $y'(\eta^+)=z'(\eta^+)$ and $y'(\zeta^-)=z'(\zeta^-)$.
Up to taking inverses of $y$ and $z$, we can assume that $y'(\eta^+)=z'(\eta^+) < 1$.
Let $\alpha$ be the first interior fixed point of $y$. Since we are looking for conjugators
fixing $\Fix(y)$ pointwise, we can restrict to find a conjugator for $y$ and $z$
in $\PL_{S,G}([\eta,\alpha])$. Since $y,z \in \PL_{S,G}^<([\eta,\alpha])$, by Proposition
\ref{thm:only-finitely-conjugators-tested} there are only finitely many candidate conjugators.
We test them and, if any of them is a conjugator in $\PL_+([\eta,\alpha])$, we extend it to $J$ through
the Stair Algorithm and test it on $J$. By the straightforward
analogues of Lemma~\ref{thm:bound-slopes-PL_+} and Proposition
\ref{thm:only-finitely-conjugators-tested} for $\PL_{S,G}^0(J)$, we can enumerate all possible conjugators.
\end{proof}

\theoremx{
\label{solofcongproblem}
The group $\PL_{S,G}(J)$ has solvable conjugacy problem. Moreover we can construct and enumerate all possible conjugators.
}

\begin{proof}
We use Theorem~\ref{congdiagonals} and suppose that
$\partial_S \Fix(y) = \partial_S \Fix(z) = \{ \eta = \alpha_0  < \alpha_1
< \ldots < \alpha_r < \alpha_{r+1} = \zeta\}$.
Now we restrict to an interval $[\alpha_i,\alpha_{i+1}]$ and consider $y,z \in \PL_{S,G}^0([\alpha_i,\alpha_{i+1}])$.
If $\Fix(y)$ contains a subinterval of $[\alpha_i,\alpha_{i+1}]$, then we must have $y=z=id$ on
the whole interval $[\alpha_i,\alpha_{i+1}]$
and so any function $g \in \PL_{S,G}([\alpha_i,\alpha_{i+1}])$ will be a conjugator. Otherwise,
$\Fix(y)$ does not contain any subinterval of $[\alpha_i,\alpha_{i+1}]$ and so we can apply the
Lemma~\ref{thm:idea-conjugacy}.
If we find a solution on each such interval, then the conjugacy problem is solvable. Otherwise, it is not.
\end{proof}

\remarkx{
\label{thm:why-last-requirement}
For the case of Thompson's group $\PL_2(I)$ there is no need of use Lemma~\ref{thm:possible-interval-initial-slopes-PLSG},
because all possible initial slopes of $g$ must be powers of $2$.
Hence, there are only finitely many conjugators with initial slope
in $[y'(0),y'(0)^{-1}]$. We test all candidate conjugators with initial slope in $[y'(0),y'(0)^{-1}]$
to conclude the procedure.
}

The argument given to solve the conjugacy problem in $\PL_{S,G}(J)$ also works,
in much the same way, to solve the power conjugacy problem. We say
that a group $G$ has \emph{solvable power conjugacy problem} if there is an algorithm
such that, given any two elements $y,z \in G$, we can determine whether there is,
or not, a $g \in G$ and two non-zero integers $m,n$ such that $g^{-1}y^m g=z^n$,
that is, there are some powers of $y$ and $z$ that are conjugate.

\theoremx{
\label{thm:power-conjugacy}
The group $\PL_{S,G}(J)$ has solvable power conjugacy problem.
}

\begin{proof}
Again, we can use Theorem~\ref{congdiagonals},
$\partial_S \Fix(y) = \partial_S \Fix(z)$ and restrict to a smaller interval $J=[\eta,\zeta]$ with
endpoints in $S$ and such that $y,z \in \PL_{S,G}^0(J)$. If $g \in \PL_{S,G}(J)$ and $m,n$ exist then we must
have that the initial slope of $y^m$ and $z^n$ must coincide. A simple argument on
the exponent of these slopes, implies that this can happen if and only if $y^m$ and $z^n$
are both powers of a common minimal power $(y^\alpha)'(\eta) = (z^\beta)'(\eta)$. Hence
the problem can be reduced to finding whether there is a $g \in \PL_{S,G}(J)$ and an integer $k$
such that $g^{-1}y^{k \alpha} g=z^{k \beta}$. By Lemma~\ref{thm:unique-roots} (that
can be naturally generalized to $\PL_{S,G}^0(J)$; see Remark~\ref{thm:closure-inside-PLSG}), we have that this is equivalent
to finding a $g \in \PL_{S,G}(J)$ such that $g^{-1}y^{\alpha} g=z^{\beta}$. Hence
solving the power conjugacy problem is equivalent to solving the conjugacy problem
for $y^\alpha$ and $z^\beta$.
\end{proof}

\section{The $k$-Simultaneous Conjugacy Problem}
\label{sec:simconj}

We will make a sequence of reductions to solve the simultaneous conjugacy problem
in $\PL_{+}(J)$ and $\PL_{S,G}(J)$. Let $M$ denote the group $\PL_+(J)$ or  $\PL_{S,G}(J)$, which will
allow us to treat both cases together.
These reductions closely follow the solution of the ordinary conjugacy problem.
First we notice that, since we know how to solve
the ordinary conjugacy problem, then solving the $(k+1)$-simultaneous conjugacy problem is
equivalent to find a positive answer to the following problem:

\problemx{
\label{thm:equiv-prob}
Is there an algorithm such that given
$(x_1, \ldots, x_k, y)$ and $(x_1, \ldots, x_k,z)$ it can decide whether
there is a function $g \in C_{M}(x_1) \cap \ldots \cap C_{M}(x_k)$
such that $g^{-1} y g = z$?}

Since we understand the structure of the intersection of centralizers,
we are going to work on solving this last question. Our strategy now is to reduce the
problem to the ordinary conjugacy problem and to isolate a very special case
that must be dealt with.

As in the the case of the ordinary conjugacy problem the first step is to determine if
the set of fixed points can be made the same.

\lemmax{
\label{simdiagonals}
Let $x_1,\ldots,x_k,y,z \in M$. We can determine whether there is, or not,
a $g \in C=C_{M}(x_1) \cap \ldots \cap C_{M}(x_k)$ such that $g(\Fix(y))=\Fix(z)$.
}
\begin{proof}
The proof is essentially the same as that of Corollary~\ref{thm:overlap-intersection-diagonal}
on each of the intervals between two fixed points of $y$ and $z$ that are in $S$. The only new tool
required is Lemma~\ref{conginZ} on the intervals where $C$ is isomorphic to $\mathbb{Z}$. We omit
the details of this proof.
\end{proof}

\lemmax{
\label{restriced-group}
Let $x_1,\ldots,x_k,y,z \in M$. The subgroup $C'$ of elements $g$
in $C_{M}(x_1) \cap \ldots \cap C_{M}(x_k)$
such that $g(\Fix(y))=\Fix(y)$ splits as a product
$$
C' = C'_{J_1}.C'_{J_2}\dots C'_{J_k}
$$
for some disjoint intervals $J_i$ with $\cup J_i = J$
where $C'_{J_i}:= \{f \in C' \mid f(t)=t, \forall t \not \in J_i\}
= C' \cap \PL_+(J_i)$.
Moreover, each
$C'_{J_i}$ is isomorphic to either $\mathbb{Z}$, or $\PL_+(J_i)\cap M$ or
is the trivial group.
}
\begin{proof}
Similar to the proof of Proposition~\ref{thm:intersectionofcentralizers}.
\end{proof}

Using the two results we reduce the simultaneous conjugacy problem to the
case $\Fix(y)=\Fix(g)$. Again we can further reduce
to the case when both $y$ and $z$ are in $\PL_{S,G}^0(J)$, but we are restricted to
use only conjugating elements from the subgroup $C'$.
By Lemma~\ref{restriced-group} the group $C'$ splits as a product of several subgroups $C'_{J_i}$,
which lead to several cases:

\emph{Case 1. The number of intervals $J_i$ is more than $1$:}
There is an interior point $\lambda$ in $J$ which is fixed by all elements in $C'$
(since $\bigcup (\partial J_i) \not\subset \partial J$).
By Lemma~\ref{cong:unique} (which can be naturally adapted to $\PL_{S,G}(J)$; see Remark~\ref{thm:closure-inside-PLSG})
there is at most one element in $C_{\PL_{S,G}(J)}(y,z)$,
which fixes $\lambda$ and we only need to verify if this element is inside $C'$.

\emph{Case 2. The number of intervals $J_i$ is exactly $1$:}
This case breaks further into 3 subcases depending on the subgroup $C'$.

\emph{Case 2a. The group $C'$ is trivial:}
The elements $y$ and $z$ are conjugate by an element in $C'$ if and only if they are the same.

\emph{Case 2b. The group $C'$ is isomorphic to $\PL_{S,G}(J)$:}
If $C'$ is the whole group, we can simply apply the algorithm which gives the solution of
the ordinary conjugacy problem.

\emph{Case 2c. The group $C'$ is isomorphic to $\mathbb{Z}$:}
We want to see if we can solve the ordinary conjugacy problem when we have a restriction on the possible conjugators.
Let $x$ denotes the generator of $C'$, thus we want to check if there exists integer $k$
such that $x^{-k} y x^k = z$. 
By assumption both $y$ and $z$ are in $\PL_{S,G}^0(J)$, solving the
ordinary conjugacy problem we find that the set $C_{\PL_{S,G}(J)}(y,z)$ is either empty or is equal to
$$
\{ \hat{y}^i g \, \mid \, i \in \mathbb{Z}\},
$$
where $\hat{y}$ is the generator of $C_{\PL_{S,G}(J)}(y)$ and $g$ is some element which conjugates
$y$ to $z$. Thus we need to find integer solutions (or show that they do not exist) of the equation
\begin{equation}
\numberwithin{equation}{section}
\label{eq:special-case}
x^k = \hat{y}^i g.
\end{equation}
This equation can be solved using the following lemma (the proof is in Subsection~\ref{ssec:final-lemma}):
\lemmax{
\label{thm:solving*}
For any $x,\hat{y},g \in \PL_{S,G}(J)$ there is an algorithm which finds
all solutions of equation (\ref{eq:special-case}).
}

Thus in all cases we can check if there exists a conjugating element in the subgroup $C'$, which
finishes the solution of the simultaneous conjugacy problem.

The previous argument proves the following theorem:
\theoremx{
The $k$-simultaneous conjugacy problem is solvable in the group $\PL_{S,G}(J)$.
Moreover we can construct and enumerate all possible conjugators.
}

\subsection{Proof of Lemma~\ref{thm:solving*}}
\label{ssec:final-lemma}

We start by proving the Lemma for the case of $\PL_2(J)$. 
We will then explain what is required to generalize
the proof to the case of $\PL_{S,G}(J)$
\footnote{The generalization to $\PL_{S,G}(J)$ is explained in the last paragraph of the current subsection.}.
We observe that both $x$ and $\hat{y}$ are in $\PL_2^0(J)$, therefore their initial slopes are not
equal to $1$. Comparing the slopes at $\eta$ and taking logarithms we obtain
\begin{equation}
\numberwithin{equation}{section}
\label{eq:Q-linearity}
k \log_2  x'(\eta^+) = \log_2 g'(\eta^+) + i \log_2 \hat{y}'(\eta^+).
\end{equation}
This equation does not have any solution unless $\log_2 g'(\eta)$ is divisible by the
greatest common divisor of $\log_2  x'(\eta^+)$ and $\log_2 \hat{y}'(\eta^+)$. If this is the
case, an elementary number theory argument tells us that all solutions are of the form
$$
k = p_1 j + q_1 \quad \mbox{and} \quad i = p_2 j +q_2,
$$
for some integers $p_1,p_2,q_1$ and $q_2$, which reduces  equation (\ref{eq:special-case}) to
\begin{equation}
\numberwithin{equation}{section}
\label{eq:reduction}
\bar{x}^ j = \bar{y}^j \bar{g}
\end{equation}
where $\bar{x}$ and $\bar{y}$ are powers of $x$ and $\hat{y}$ respectively and
$\bar{g}'(\eta^+)=1$.

If $\Fix (\bar x) \not =\Fix (\bar y)$ we can use Lemma~\ref{conginZ}
to solve equation (\ref{eq:reduction}). We can also compare the derivatives at all fixed points and this
will give us a unique solution (or that there does not exist any solution) for $j$
unless the following are satisfied
$$
\bar{x}'(\mu)=\bar{y}'(\mu)
\quad\mbox{and}\quad
\bar{g}'(\mu),
$$
for any $\mu \in \Fix (\bar x)$. Equation (\ref{eq:reduction}) can be written as
\begin{equation}
\numberwithin{equation}{section}
\label{eq:rewriting-reduction}
 \bar{g}= \bar{x}^ j\bar{y}^{-j}.
\end{equation}
If $\bar{x} = \bar{y}$ equation (\ref{eq:rewriting-reduction}) has solutions if and only if $\bar g= id$ and in this case
any integer $j$ is a solution. Thus the only non-trivial case when $\bar{x}\not= \bar{y}$.

Without loss of generality we may assume that $\bar{x},\bar{y} \in \PL_+^<([\mu_1,\mu_2])$
for some consecutive $\mu_1$ and $\mu_2$ in $\partial \Fix (\bar x)= \partial \Fix (\bar y)$.
Let $p$ denote the function $\bar{x} \bar{y}^{-1}$ and let $\lambda$
be the closest breakpoint of $p$ to $\mu_1$, i.e., $p(t) = t$ for all $\mu_1 \leq t \leq \lambda$ and
$p(\lambda + \varepsilon) \not = \lambda + \varepsilon$ if $\varepsilon>0$ is sufficiently small.
For any $j > 0$ we can write
\begin{equation}
\numberwithin{equation}{section}
\label{eq:-second-rewriting}
\bar{g}= \bar{x}^ j\bar{y}^{-j} = p p^{\bar y^{-1}} \dots p^{\bar y^{-j+1}}.
\end{equation}
It is clear that the first breakpoint for $p^{\bar y^r}$, for any integer $r$,
is given by $\bar{y}^{-r}(\lambda)$. Since $\bar{y} \in \PL_+^<([\mu_1,\mu_2])$,
formula (\ref{eq:-second-rewriting}) gives that the first breakpoint of
$\bar{g}$ is at $\bar{y}^{j-1}(\lambda)$. There can be at most one positive $j$ such that the number
$\bar{y}^{j-1}(\lambda)$ coincides with the actual first breakpoint of $\bar{g}$.
Therefore, we can find if equation (\ref{eq:reduction}) has solutions for
positive  $j$.
If $j$ is negative we can similarly write,
\begin{equation}
\numberwithin{equation}{section}
\label{eq:last-rewriting}
\bar{g}^{-1}= \bar{y}^{-j}\bar{x}^{j} = \bar{p} \bar{p}^{\bar x^{-1}} \dots \bar{p}^{\bar x^{j+1}}.
\end{equation}
where $\bar{p} := p^{-1}$. Since $\bar{x} \in \PL_+^<([\mu_1,\mu_2])$,
formula (\ref{eq:last-rewriting}) gives that the first breakpoint of
$\bar{g}^{-1}$ is at $\bar{x}^{-j-1}(\lambda)$. Therefore, we can find if equation (\ref{eq:reduction}) has solutions for
negative  $j$.

This completes the proof of Lemma~\ref{thm:solving*} for $\PL_2(J)$.
To generalize this proof to the groups $\PL_{S,G}(J)$ we observe
that all of the previous proof has been carried out in $\PL_+(J)$, save for the first step, that is taking logarithms
to get an argument to pass from 
equation~(\ref{eq:special-case}) to equation (\ref{eq:reduction}). 
To do this step in $\PL_{S,G}(J)$,
we appeal to the last of the requirement in section~\ref{sec:comp_req}.

\section{Interesting Examples}
 \label{sec:examples}

Now that we have developed the general theory, we are going to see a few interesting examples where
the simultaneous conjugacy problem is solvable.
We will not dwell too much on the details here, sketching only why it is possible to
verify the requirements.

\examplex{\label{thm:example-1}
$S=\mathbb{Q}$ and $G=\mathbb{Q}^*_{\ge 0}=\mathbb{Q} \cap (0,\infty)$.

There are many structures which can be used to represent the rational numbers,
which comes with algorithms for performing the arithmetic operations, which give
us the oracles in the first group. The oracles in the second group are very easy 
to implement since $\mathbb{Q}$ is a field and the quotients 
$S / \mathcal{I}=\{0\}$ and $S / (t-1)\mathcal{I}=\{0\}$
consist of just one element.
The last oracle which is need for solving the simultaneous conjugacy problem is 
slightly more complicated -- we need to factor $a,b,c$ as product of prime 
numbers and then reduce the problem to solving several congruences in integers. 
}

\examplex{\label{thm:example-2}
$S$ finite real algebraic extension over $\mathbb{Q}$
and $G=S^*:= S \cap (0,\infty)$. This is the same as the previous example, we
only need to ``implement'' the field $S$.
}

\examplex{
\label{thm:example-rationals}
$S=\mathbb{Z}\big[\frac{1}{n_1}, \ldots, \frac{1}{n_k}\big]$
and $G=\langle n_1, \ldots, n_k \rangle$ for $n_1,\ldots,n_k \in \mathbb{Z}$.

As in Example \ref{thm:example-1} there are many data structures to represent
$S$ and $G$, which provide the oracles in the first group.
For the oracles in the second group one observes that 
$S / \mathcal{I} \cong \mathbb{Z}/d\mathbb{Z}$, where
$d:=\mathrm{GCD}(n_1-1, \ldots, n_k-1)$.
This reduces an effective solution of the membership problem in 
$\mathcal{I}$ to expressing a given element in
$d\mathbb{Z}$ as sum of multiples of $n_i-1$, which can be done
using the Euclid's algorithm.
As in the previous example the implementing the last oracle relies on
the factorization of of integers as product of primes.
For $k=1$, we recall that the groups $\PL_{S,G}(I)$ are known as 
\emph{generalized Thompson's groups}.
}
\examplex{$S=\mathbb{Z}\big[\frac{1}{n_1}, \ldots, \frac{1}{n_k}, \ldots \big]$
with $G=\langle \{n_i\}_{i \in \mathbb{N}} \rangle$, 
where $\{n_i\}_{i \in \mathbb{N}} \subseteq \mathbb{Z}$.

This example can be reduced to the previous one. If we are given a finite set
$E$ of elements in $\PL_{S,G}(I)$ we can consider the set $\{n_{i_1}^{\alpha_{i_1}},\ldots,
n_{i_v}^{\alpha_{i_v}}\}$ of all slopes of elements of $E$. Then $E \subseteq \PL_{S',G'}(I)$
where $S':=\mathbb{Z}\big[\frac{1}{n_{i_1}}, \ldots, \frac{1}{n_{i_v}}\big]$ and
$G':=\langle n_{i_1}, \ldots, n_{i_v} \rangle$. By Corollary \ref{thm:conjugator-still-in-PLSG}
we know that if there is a conjugator, it must be in $\PL_{S',G'}(I)$.
}

\remarkx{
\label{thm:multiple-orbits-of-points}
In general, given two intervals $J_1,J_2$ with endpoints in $S$, it is not clear whether or not the groups
$\PL_{S,G}(J_1)$ and $\PL_{S,G}(J_2)$ are isomorphic.
Proposition~\ref{thm:equiv-cond-trans} tells us that two elements in $S$ are in the same $\PL_{S,G}$-orbit if their image
under the map $\pi$ is the same.
For example in the cases
$S=\mathbb{R},G=\mathbb{R}_+$ and $S=\mathbb{Q},G=\mathbb{Q}^*$ and $S=\mathbb{Z}\left[\frac{1}{2}\right],G=\langle 2 \rangle$,
it is not difficult to see that every two points in $S$ have the same image under $\pi$ (the case of $F$ is treated in 
Lemma \ref{thm:standard-cannon-floyd-parry})
and that any two groups $\PL_{S,G}(J_1)$ and $\PL_{S,G}(J_2)$ are thus isomorphic, for any two intervals
$J_1,J_2$ with endpoints in $S$. In fact, if there is a $\PL_{S,G}$-map $\varphi:J_1 \to J_2$, then conjugation by $\varphi$
yields an isomorphism between $\PL_{S,G}(J_1)$ and $\PL_{S,G}(J_2)$.

On the other hand, if we consider generalized Thompson's groups
(see Example~\ref{thm:example-rationals}) and use the map $\pi$, 
it is straightforward to show that the number of orbits of elements is finite but more than one, for certain
choices of $n_1,\ldots,n_k$
(see Example~\ref{thm:example-rationals} for a proof of this),
hence there are only finitely many inequivalent types intervals $J$ with endpoints in $S$.
This implies that there can be at most only finitely many isomorphism classes for the groups $\PL_{S,G}(J)$,
for $S=\mathbb{Z}\big[\frac{1}{n_1}, \ldots, \frac{1}{n_k}\big]$
and $G=\langle n_1, \ldots, n_k \rangle$ for $n_1,\ldots,n_k \in \mathbb{Z}$.
We observe that the generalized Thompson's groups which are most often studied are those where we assume
that $\mathrm{GCD}(n_1-1,\ldots,n_k-1)=1$, which implies that $S/\mathcal{I}$ is trivial.
In general, it seems likely that if two elements $\alpha,\beta \in S$ have different image under $\pi$
then the groups $\PL_{S,G}([0,\alpha])$ and $\PL_{S,G}([0,\beta])$ are not isomorphic, but it is not easy to prove it.
}

\section*{Acknowledgments}

The first author was partially supported by AMS Centennial Fellowship.
This work is part of the second author's PhD thesis at Cornell University.
The second author gratefully acknowledges the Centre de Recerca Matem\`atica (CRM)
and its staff for the support received during the completion of this work.
The authors would like to thank
Collin Bleak, Martin Bridson, Ken Brown, Kai-Uwe Bux and
Mark Sapir for the useful discussions and
comments on earlier drafts of this paper.


\end{document}